\documentclass[12pt]{scrartcl}
\usepackage{a4}
\usepackage{amsthm}
\usepackage{amsmath}
\usepackage{amssymb}
\usepackage{amsfonts}
\usepackage{mathrsfs}
\usepackage{dsfont}
\usepackage{latexsym}
\usepackage{color}
\usepackage{longtable}
\usepackage{bbm,exscale}
\usepackage[titletoc]{appendix}
\definecolor{Myblue}{rgb}{0,0,0.6}
\usepackage[a4paper,colorlinks,citecolor=Myblue,linkcolor=Myblue,urlcolor=Myblue,pdfpagemode=None]{hyperref}

\usepackage[a4paper,margin=2.5cm]{geometry}

  \vfuzz \hfuzz
\newcommand{\be}{\begin{equation}}
\newcommand{\ee}{\end{equation}}
\newcommand{\bes}{\begin{equation*}}
\newcommand{\ees}{\end{equation*}}
\newcommand{\bea}{\begin{eqnarray}}
\newcommand{\eea}{\end{eqnarray}}
\newcommand{\ba}{\begin{array}}
\newcommand{\ea}{\end{array}}

\def\la {\label}

\newcommand{\C}{\mathbb{C}} 

\newcommand{\Q}{\mathbb{Q}}

\newcommand{\Z}{\mathbb{Z}} 

\def\1{\ifmmode\mathrm{1\!l}\else\mbox{\(\mathrm{1\!l}\)}\fi}

\newcommand{\oeq} {\sim_\mathrm{oeq}} 

\def\maxm{{\mathfrak m}}

\allowdisplaybreaks

\deffootnote[1em]{1em}{1em}{\textsuperscript{\thefootnotemark}}

\numberwithin{equation}{section}

\begin{document}

\title{Orbifold equivalence: structure and new examples}

\author{Andreas Recknagel$$ \quad Paul Weinreb
\\
\\[-0.1cm] 
 {\normalsize \textsl{King's College London, Department of Mathematics,}} \\ [-0.1cm]
   {\normalsize \textsl{Strand, London WC2R\,2LS, UK}}\\[-0.1cm]
}

\date{}
\maketitle

\footnote{\texttt{andreas.recknagel@kcl.ac.uk, paul.weinreb@kcl.ac.uk}}

\begin{abstract}
  \noindent 
 Orbifold equivalence is a notion of symmetry that does not rely
 on group actions. Among other applications, it leads to surprising
 connections between hitherto unrelated singularities.
 While the concept can be defined in a very general category-theoretic
 language, we focus on the most explicit setting in terms of
 matrix factorisations, where orbifold equivalences arise from defects 
 with special properties. Examples are relatively difficult to construct,
 but we uncover some structural features that distinguish orbifold
 equivalences -- most notably a finite perturbation expansion. We use
 those properties to devise a search algorithm, then present
 some new examples including Arnold singularities.

 
\end{abstract}

\newpage

\tableofcontents

\vspace{1cm}

\section{Introduction}\label{sec:introduction}

Orbifold equivalence is a phenomenon discovered when trying to describe 
well-known ideas connected with the action of symmetry groups (originally
in quantum field theories) in terms of abstract category-theoretic terms.
It turned out that all the data one is interested in when studying these
``orbifolds'' can be extracted from a separable symmetric
Frobenius algebra -- which may, but need not, arise from a group action. 
The abstraction therefore provides a generalised notion of symmetry which
does not rely on groups.

The original formulation led to some rather strong results concerning the
classification of rational conformal field theories \cite{FFRS}, the
more abstract bicategory point of view taken in \cite{CR3} allows for
a wider range of applications, including so-called topological
Landau-Ginzburg models.

In the latter context, orbifold equivalence provides a novel equivalence
relation for quasi-homogeneous polynomials,
leading to unexpected relations between singularities -- e.g.\ between
simple singularities of types $A$ and $E$ -- and to equivalences of
categories  associated with them, such as categories of matrix factorisations
and of representations of quiver path algebras. The new equivalence also 
implies ``dualities'' between different topological field theories (correlation
functions of  one model can be computed in another).

\vspace{.1cm} 
The main aim of the present paper is to construct explicit examples of
orbifold equivalences, which in the Landau-Ginzburg context are nothing but 
matrix factorisations with special properties.

\vspace{.3cm} 
We start by recalling some basic definitions concerning matrix factorisations,
then make a few remarks about the mathematics and physics context. The
``special property'' we demand (namely invertible quantum dimensions) will be
addressed in section 2. Our main new results will be presented in sections 3
and 4: We uncover some structural properties of orbifold equivalences (they
have a finite graded perturbation expansion), exploit this to set up
a search algorithm, then list examples of orbifold equivalences found in this way. 

\vspace{.3cm} 
A rank $N$ \textsl{matrix factorisation} of a polynomial
$W \in 
{\mathbb C}[z_1,\ldots,z_k]$
is a pair of $N \times N$ matrices $E,\;J$ with polynomial entries satisfying
\bea
  E\,J = J\,E = W\, \mathbf{1}_N\ .
\eea
We can collect $E$ and $J$ into a single matrix $Q$,  
\bea \la{QEJ}
   Q= \left(\begin{matrix}\ 0&E\ \\  J&0\  \end{matrix} \right)      
   \in M_{2N}({\mathbb C}[z])\ ,
\eea
satisfying $Q^2 = W\, \mathbf{1}_{2N}$; we have abbreviated $z = (z_1,\ldots,z_k)$.  
This notation is one way to make the inherent $\Z_2$-grading of matrix factorisations
explicit: $Q$ anti-commutes with 
 \bea \la{sigmadef}
 \sigma = \left(\begin{matrix}\  \mathbf{1}_{{N}}  & 0\ \\  0 &  -\mathbf{1}_{{N}}\ 
  \end{matrix} \right)
\eea
and we can use $\sigma$ to decompose the space $M_{2N}({\mathbb C}[z])$ into
even and odd elements: The former commute with $\sigma$ and are called
``bosonic'' in the physics context, the latter anti-commute with $\sigma$
and are referred to as ``fermionic''. 

Every polynomial $W$ admits matrix factorisations: Any
factorisation of $W(z_1) = z_1^p$ provides a rank 1 factorisation, 
and once we have matrix factorisations $Q_a$ of $W_a$ and $Q_b$ of $W_b$,
the so-called tensor product factorisation
$Q = Q_a\,\hat\otimes\,Q_b$ provides a  matrix factorisation of
$W = W_a + W_b$.  This $Q$ is formed as in (\ref{QEJ}) from 
$$
E :=
  \left( \begin{array}{cc}
         J_a \otimes\mathbf{1}  & -\mathbf{1}\otimes E_b \\
    \mathbf{1}\otimes J_b  &    E_a \otimes\mathbf{1}   
  \end{array} 
  \right)\ , \quad\quad   
  J :=
  \left( \begin{array}{cc}
      E_a \otimes\mathbf{1}  & \mathbf{1}\otimes E_b \\
    -\mathbf{1}\otimes J_b  &    J_a \otimes \mathbf{1}   
  \end{array} 
 \right)         \ .
$$
Unfortunately, it will turn out (see section 3) that matrix factorisations
obtained as iterated tensor products of factorisations of monomials of $W$
in general do not have the special additional properties we are interested in. 

\vspace{.2cm} 
In the following, we will exclusively focus on \textsl{graded} matrix factorisations:

First of all, we assume that the polynomial $W(z)$ is quasi-homogeneous, i.e.\
that there exist rational numbers $|z_i| > 0$, called the weights of $z_i$,  
such that for any $\lambda \in \C^{\times}:=\C \setminus \{0\}$ we have
$$
W(\lambda^{|z_1|} z_1,...,\lambda^{|z_k|} z_k)\ =\ \lambda^{D_W} W(z_1,...,z_k)
$$
for some $D_W \in \Q_+$, the weight of $W$. Unless specified otherwise, we will assume that $D_W = 2$, and also that $W \in \maxm^2$ where $\maxm = \langle z_1,\ldots,z_k\rangle$
is the maximal ideal of $\C[z]$. For some applications, it is important that
the Jacobi ring
$\mathrm{Jac}(W) = \C[z]/\langle \partial_{z_1}W,\ldots,\partial_{z_k}W \rangle$
is finite-dimensional as a $\C$-vector space, so let us assume this.

We will refer to quasi-homogeneous polynomials of weight 2 as ``potentials'',
and to the rational number 
\bea\la{cench} 
\hat c(W) := \sum_{i=1}^k(1-|z_i|) 
\eea
as the \textsl{central charge} of the potential $W$.  

In the following, let us use the abbreviation
\be\la{lambdaaction}
\lambda \rhd z := (\lambda^{|z_1|} z_1,...,\lambda^{|z_k|} z_k)
\ee
for the $\C^\times$-action.

We call a rank $N$ matrix factorisation $Q$ of a potential $W$ \textsl{graded} if
there exists a diagonal matrix (the ``grading matrix'' of $Q$) 
$$U(\lambda)\ =\ \mathrm{diag}(\lambda^{g_1},\ldots,\lambda^{g_{2N}})
$$
with $g_i \in \Q$ such that 
\be\la{grading}
U(\lambda)\,Q(\lambda \rhd z )\,U(\lambda)^{-1} = \lambda\; Q(z)
\ee
for all $\lambda\in\C^\times$. We can set $g_1=0$ without loss of generality.

For simplicity, we will assume that $Q(z)$ has no non-zero constant entries. 
Otherwise $Q$ is decomposable: using row and column transformations, it
can be brought into the form
$\widetilde{Q} \oplus Q_\mathrm{triv}$ where $Q_\mathrm{triv}$ is the trivial rank 1
factorisation $W = 1 \cdot W$.

\vspace{0.3cm} 
Probably the first, very simple, example of a matrix factorisation made
its appearance in the Dirac equation, but only with Eisenbud's discovery
that free resolutions of modules over $\C[z]/\langle W \rangle$ become
periodic \cite{Eis} was it realised that matrix factorisations can be defined
for general potentials and are a useful tools in mathematics.

Later, a category-theoretic point of view was introduced: one
can e.g.\ form a category $\mathrm{hmf}^\mathrm{gr}(W)$ whose objects are (finite
rank) graded matrix factorisations of $W$, and where (even or odd) morphisms
are given by the cohomology of the differential
$d_{Q_1Q_2}$ acting as    
\be  \la{DQops}
d_{Q_1Q_2}(A) = Q_1A- (-1)^{s(A)}\,AQ_2
\ee
on $A\in M_{2N}(\C[z])$, where $s(A)=0$ if $A$ is even with respect to the
$\Z_2$-grading $\sigma$ in (\ref{sigmadef}), and $s(A)=1$ if $A$ is odd.

It was shown that $\mathrm{hmf}^\mathrm{gr}(W)$ is equivalent to the derived
category of coherent sheaves on zero locus of $W$ -- and also to categories of
maximal Cohen-Macaulay modules, and of quiver path algebras; see in particular
\cite{Orlov}, but also \cite{Buchw}. Among these categories, $\mathrm{hmf}^\mathrm{gr}(W)$
is the one where explicit computations are easiest to perform.

\vspace{0.2cm}
Let us also make some remarks on the most notable application of matrix factorisations
in physics, namely topological Landau-Ginzburg models.  (We include these
comments mainly because this context is the origin of some of the terminology;
an understanding of the physical concepts is not required for the remainder
of the paper.)
Topological Landau-Ginzburg models are  supersymmetric quantum field theories
on a two-dimensional worldsheet; $W(z)$ appears as interaction potential,
the degrees of freedom in the interior of the world-sheet
(the ``bulk fields'') are described by the Jacobi ring
$\mathrm{Jac}(W)$. 

Within the physics literature, there is strong evidence for a relation
between  supersymmetric Landau-Ginzburg models with potential $W$ and
supersymmetric conformal field theories where a Virasoro algebra with
central charge $c = 3 \,\hat c(W)$ acts. The conformal field theory 
is thought to describe the ``IR renormalisation group fixed point''
of the Landau-Ginzburg model and motivates the term ``central charge''  for
the quantity (\ref{cench}). 

\vspace{.1cm} 
Additional structure appears if the worldsheet of the Landau-Ginzburg
has boundaries: the possible supersymmetry-preserving boundary conditions
are precisely
the matrix factorisations $Q$ of $W$ \cite{KL0,BHLS,KL2,HerLaz}, and bosonic 
resp.\ fermionic degrees of freedom on the boundary are given by the
even resp.\ odd cohomology $H^\bullet_Q$ of the differential $d_{QQ}$ defined
in eq.\ (\ref{DQops}). 

Correlation functions in topological Landau-Ginzburg are computed
as residues of functions of several complex variables (see
e.g.\ \cite{GriffithsHarris} for many details): If the
worldsheet has no boundary, the correlator of any element
$\phi \in \mathrm{Jac}(W)$ is \cite{Vafa} 
\be\la{bulkcorr}
\langle\,\phi\,\rangle_W =
       \mathrm{res}_z \biggl \lbrack
       { \phi  \over \partial_{z_1}W \cdots \partial_{z_k}W } \biggr \rbrack
       \ .
\ee
In a model where the worldsheet has a boundary with boundary
condition described by a matrix factorisation
$Q$ of $W$, one has \cite{KL1,HerLaz}
\be\la{bdycorr}
\langle\;\phi\,\psi\;\rangle_Q^\mathrm{KapLi} =
     \mathrm{res}_z\biggl\lbrack
  { \phi\; \mathrm{str} \bigl( \partial_{z_1}Q \cdots \partial_{z_k}Q \; \psi \bigr)
    \over \partial_{z_1}W \cdots \partial_{z_k}W } \biggr\rbrack   
\ee
for any bulk field $\phi \in \mathrm{Jac}(W)$ and any boundary field
$\psi \in H^\bullet_Q$. The supertrace is defined using the
$\Z_2$-grading from (\ref{sigmadef}), as
$\mathrm{str}(A) := \mathrm{tr}(\sigma\,A)$.

The formula (\ref{bdycorr}) is often referred to
as Kapustin-Li correlator; a closely related expression will be used,
in the next section, to define the ``special property'' the
matrix factorisations of our interest are required to have. 

The correlations functions above were first computed in physics,
via localisation of path integrals for supersymmetric topological
quantum field theories, but they have since been discussed in
purely mathematical terms, notably in \cite{MurfetResidues,DyckMuKapLi}.

\vspace{.2cm} 
Instead of worldsheets with boundary, one can also consider worldsheets
which are divided into two domains by a ``fault line'', and the degrees
of freedom on the two sides may be governed by two different Landau-Ginzburg
potentials $V_1(x)$ and $V_2(y)$. Such an arrangement is called a (topological)
\textsl{defect}, and is described \cite{BR} by a matrix factorisation $Q(x,y)$ of
$V_1(x)-V_2(y)$. Degrees of freedom localised on the defect line are described
by the morphisms (bosonic or fermionic) of $Q(x,y)$, analogously to the
boundary case. Boundary conditions of a Landau-Ginzburg model with
potential $V_1(x)$ can be viewed as defects between $V_1(x)$ and the
trivial model $V_2=0$. 

\vspace{.2cm}
Topological defects come with additional structure, called the
fusion product:
In the fault line picture, two defect lines which partition 
a worldsheet into three regions, with
potentials $V_1(x)$, $V_3(x')$ in the outer regions and $V_2(y)$ in
the middle, can be moved on top of each other, leaving a single defect
between $V_1(x)$ and $V_3(x')$. In terms of matrix factorisations,
the tensor product $Q_{12}(x,y)\,\hat\otimes\,Q_{23}(y,x')$ of two
matrix factorisations $Q_{12}(x,y)$ of $V_1(x)-V_2(y)$ and $Q_{23}(y,x')$
of $V_2(y)-V_3(x')$ is a matrix factorisation of $V_1(x)-V_3(x')$. 
This has infinite rank over $\C[x,x']$, but is equivalent (by a similarity
transformation) to a finite-rank defect \cite{BR} depending on $x,x'$ only;
extracting this finite rank defect yields a representative of the fusion
product $Q_{12}\star Q_{23}$. The full construction is somewhat technical
(it involves finding and splitting an idempotent morphism of the tensor
product, see \cite{CM1} for details), but implementable on a computer.
The construction described in \cite{MurfetCuts} should prove easier
to apply. 

A concrete mathematical application of defects appeared in the work of
Khovanov and Rozansky \cite{KhR}, who proposed to use matrix
factorisations to define link invariants that generalise
those of Reshetikhin-Turaev (``categorification of the Jones polynomial'').
These invariants were made explicitly computable using the fusion product
in \cite{CM1}. 

\vspace{.2cm}
At a more abstract level, topological defects in Landau-Ginzburg models,
together with structures such as the fusion product, 
form a bicategory ${\cal LG}$, where objects are given by Landau-Ginzburg
potentials, 1-morphisms by defects between two potentials (i.e.\ matrix 
factorisations of the difference),
and 2-morphisms by morphisms (as occurring in (\ref{DQops})) of those matrix
factorisations. This bicategory is ``graded pivotal'' \cite{CM2,CR3,Cr},
in particular it has adjoints: for each 1-morphism $Q$, i.e.\ each defect
between $V_1(x)$ and $V_2(y)$, the adjoint $Q^\dagger$, a defect between
$V_2(y)$ and $V_1(x)$, is given by
\be\la{Qadj}
Q^\dagger = \left(\begin{matrix}\ 0&J^T\ \\\ -E^T&0\  \end{matrix} \right)\ .         
\ee
(This equation holds if the number of $y$-variables is even, otherwise
there is an additional exchange of $E$ and $J$; see \cite{CM2,CR3,CM3} for
details, which will not play a role in what follows.) 

\vspace{.2cm}
A detailed knowledge of category theory is not required to understand the
results of the present paper. Indeed, while the category framework is
convenient, perhaps even indispensable, to develop the notion of orbifold
equivalence and to fully appreciate its wide-ranging applications (including
possible extension to higher-dimensional topological field theories), the
pedestrian approach via explicit matrix factorisations seems much better suited
to construct examples.

\section{Orbifold equivalence}\label{sec:OEQ}
\subsection{Definition and general properties}\label{OEQ2}

We now come to the definition of the ``special property'' we require the
defects of interest to have. In the following, let $V_1 \in \C[x]$ and
$V_2 \in \C[y]$,  where $x = (x_1,\ldots,x_n)$ and
$y = (y_1,\ldots,y_m)$, be two potentials.

\vspace{0.3cm}\noindent
 \textbf{Definition 2.1: }
$V_1$ and $V_2$ are \textsl{orbifold equivalent} if there 
exists a (graded) matrix factorisation $Q(x,y)$ of $V_1(x)-V_2(y)$ for which
the quantum dimensions $q_L(Q)$ and $q_R(Q)$ are invertible.

\noindent
The \textsl{quantum dimensions} of $Q$ are defined as
\bea \la{qdimdef}
q_L(Q)&=&(-1)^{{m+1 \choose 2}}\
\mathrm{res}_x \Big{[}\,\frac{\mathrm{str}(\partial_{x_1}Q\cdots\partial_{x_n}Q\cdot
    \partial_{y_1}Q\cdots\partial_{y_m}Q)}
  {\partial_{x_1}V_1 \cdots \partial_{x_n}V_1}\Big]
\cr\cr\cr 
q_R(Q)&=&(-1)^{{n+1 \choose 2}}\ 
\mathrm{res}_y \Big{[}\,\frac{\mathrm{str}(\partial_{x_1}Q\cdots\partial_{x_n}Q\cdot
    \partial_{y_1}Q\cdots\partial_{y_m}Q)}
  {\partial_{y_1}V_2 \cdots \partial_{y_m}V_2}\Big]
\eea
We have used $\mathrm{str}(A) = \mathrm{tr}(\sigma\,A)$ to abbreviate the
supertrace, defined with the help of the $\Z_2$-grading from
(\ref{sigmadef}). 

\noindent
We will call such a $Q$ an orbifold equivalence between $V_1$ and $V_2$,
and we will write $V_1 \oeq V_2$ to indicate that $V_1$ and $V_2$ are
orbifold-equivalent. 
\vspace{.3cm} 

We will see below that, in the graded case, the quantum dimensions are
complex numbers, so ``invertible'' simply means ``non-zero''. 

One can show that the quantum dimensions are invariant under permutations of
the variables up to a sign \cite{DyckMuPushing}. Notice also the close
similarity of (\ref{qdimdef}) to the Kapustin-Li correlator; this will
play a role in section 3.

One could in principle drop the requirement that $Q$ is a graded matrix
factorisation (or even that $V_1$ and $V_2$ are quasi-homogeneous), the
quantum dimensions can be computed for any matrix with polynomial entries. 
Not much is known in this general situation, so we
restrict ourselves to quasi-homogeneous potentials and
graded matrix factorisations in this paper.

\vspace{.2cm} 

We summarise some abstract properties of the notions of orbifold equivalence 
and quantum dimensions in a theorem; all statements were proven before,
see \cite{CM2,CR3,CRCR} and references therein: 

\vspace{0.2cm}\noindent
\textbf{Theorem 2.2: }
\begin{itemize} 

\item[(a)] $\oeq$ is an equivalence relation.

\noindent
\item[(b)] $V_1 \oeq V_1 + y_1^2 + y_2^2$ (Kn\"orrer periodicity).

\noindent
\item[(c)] If $V_1(x) \oeq V_2(y)$ and $V_3(x')\oeq V_4(y')$, then
$V_1(x)+V_3(x') \oeq V_2(y)+V_4(y')$.
(Note that in this relation each potential depends on a
separate set of variables.)

\noindent
\item[(d)] The quantum dimensions do not change under similarity transformations,
i.e.\ $q_L(Q) = q_L(U\,Q\,U^{-1})$ for any invertible even matrix $U$;
analogously for $q_R(Q)$.

\noindent
\item[(e)] The quantum dimensions are additive with respect to forming direct
  sums: if
$Q$ and $\widetilde{Q}$ are two matrix factorisations of $V_1 - V_2$, then 
$q_{L}(Q \oplus \widetilde{Q}) = q_{L}(Q) + q_{L}(\widetilde{Q})$, and analogously
for $q_R(Q \oplus \widetilde{Q})$.

\noindent
\item[(f)] Up to signs, the quantum dimensions are multiplicative with
  respect to fusion products $Q \star \widetilde{Q}$, and with
  respect to forming tensor products $Q_{12}(x,y)\,\hat\otimes\,Q_{34}(x',y')$
  (where $Q_{12}$ factorises $V_1(x)-V_2(y)$ and $Q_{34}$ factorises
  $V_3(x')-V_4(y')$, cf. item c).

\noindent
\item[(g)] Passing to the adjoint defect interchanges left and right
  quantum dimensions: $q_L(Q^\dagger) = q_R(Q)$ and $q_R(Q^\dagger) = q_L(Q)$.

\end{itemize}

\vspace{.3cm} 
A quantity of central importance in the bicategory treatment of orbifold
equivalences is $A(Q) := Q^\dagger \star Q$, sometimes
called ``symmetry defect''. This fusion product is a defect from
$V_2$ to itself, and it can be shown \cite{CR3} that
\be\la{hmfmodequiv}   
\mathrm{hmf}^\mathrm{gr}(V_1) \simeq \mathrm{mod}(Q^\dagger\star Q)
\ee
where the right hand side denotes the category of modules over $A(Q)$,
consisting of matrix factorisations of $V_2$ on which the defect
$A(Q)$ acts via the fusion product. This equivalence of categories
is one of several relations existing between structures associated to $V_1$
and to $V_2$ as soon as the two potentials are orbifold equivalent.

Within the domain of Landau-Ginzburg models, orbifold equivalence leads
to a ``duality'' of the two topological field theories: bulk correlators
in the $V_1$-model can be computed as correlators in the $V_2$-model enriched
by defect lines (the defect being $A(Q)$ -- see e.g.\ \cite{CR3} for a nice
pictorial presentation of this fact.

All one needs to prove these statements in the bicategory language is that
$A(Q)$ is a ``separable symmetric Frobenius algebra''. It is this power
of abstraction that made it possible to realise that features one is
familiar with from orbifold groups can persist without groups being involved.

\vspace{.2cm}
Explicit computations involving $A(Q)$
can become rather tedious when dealing with complicated orbifold equivalences. 
However, there is a very simple numerical \textsl{invariant} which contains
useful information, namely the (left or right) quantum dimension of $A(Q)$:
Using the facts collected in Theorem 2.2, we find 
\be\la{qlr}
q_L(A(Q)) = q_L(Q^\dagger) q_L(Q) = q_L(Q)q_R(Q) = q_R(A(Q))\ .
\ee
E.g., if $Q$ is an indecomposable defect with $q_L(A(Q)) \neq \pm 1$, then
$Q$ is a ``true orbifold equivalence'' rather than a ``mere equivalence''
in the bicategory ${\cal LG}$, i.e.\ there cannot be a $\widetilde{Q}$ such
that $Q\star \widetilde{Q}$ and $\widetilde{Q}\star Q$ are
similar to the unit 1-morphisms $I_{V_1}$ and $I_{V_2}$ of $V_1$ resp.\ $V_2$ (see the
next subsection for the definition of $I_V$).

Perhaps more interestingly, an orbifold equivalence $Q$ does not arise from the
action of a finite symmetry group on the potential unless $q_L(A(Q))$ is contained
in some cyclotomic field: this follows from constraints on the quantum dimensions
of orbifold equivalences associated with group actions, see the remarks in the
next subsection.

\vspace{.2cm} 

By definition, orbifold equivalence describes a property of a pair of
potentials, a defect between them with non-vanishing quantum dimensions
merely needs to exist. Ultimately, one would like to be able to read off
directly from the potentials whether they are orbifold equivalent or not.
So far, however, only the following two facts are known to be necessary
criteria for orbifold equivalence:

\vspace{0.2cm}\noindent
\textbf{Proposition 2.3: } Using the notations from Def.\ 2.1, $V_1(x)$ and
$V_2(y)$ are orbifold equivalent only if the total number of variables
$n+m$ is even and only if the two potentials have the same central charge, 
$\hat c(V_1) = \hat c(V_2)$.

\vspace{.2cm}\noindent
Both statements were proven in \cite{CR3}, we will give a slightly modified
derivation in section 3.
%
%
%
The first criterion of Prop.\ 2.3 
is easy to circumvent: if $V_1$ and $V_2$ have an odd total number of variables, one
can pass from $V_1$ to $V_1+x_{n+1}^2$, which has the same central charge.
The second criterion, on the other hand, is definitely necessary for the
existence of (graded) orbifold equivalences, and a natural question is whether
it is already sufficient. We will make a few comments on this in section 5.

\subsection{Known examples}\label{sec:examples}

We briefly recapitulate the examples of orbifold equivalences known so far.

\vspace{0.1cm}
Reflexivity of $\oeq$ is ensured by the identity defect%
, an orbifold equivalence $I_V(x,y)$ between $V(x)$ and $V(y)$ with
quantum dimensions $q_L(I_V) = q_R(I_V) = 1$. This takes the form of a nested
tensor product $I_V:=Q_{(1)}\,\hat\otimes\cdots\hat\otimes\, Q_{(n)}$ of rank 1
matrix factorisations with
$$
E_{(i)}=[V(x_1,...,x_i,y_{i+1},...y_n)-V(x_1,...,x_{i-1},y_{i},...,y_n)]/(x_i-y_i)\ ,
\ \ J_{(i)}=x_i-y_i\ .
$$
In the bicategory picture, the identity defect provides the unit (with
respect to the fusion product) 1-morphism in $\mathrm{End}(V)$. In the Landau-Ginzburg
context, $I_V$ can be thought of as an ``invisible defect''. 

The example from which ``orbifold equivalences'' derive their name involves
a symmetry group $G$ of the potential $V$, i.e.\ a finite subgroup of
$\C[x]$-automorphisms which leaves $V$ invariant. Then one can, for
each $g\in G$, construct ``twisted'' identity defects  $I^g_V$ formed like
$I_V$ above, with $J_i$ replaced by $J^g_i =x_i-g(y_i)$, and $E_i$
replaced accordingly. Details are given in  \cite{BR,CR3,BCP1,BCP2},
where it is also shown that the quantum dimensions of $I^g_V$ are
given by $\det(g)^{\pm1}$ (hence contained in the cyclotomic field
determined by the order of the symmetry group $G$). 
In this special situation, the symmetry defect $A(Q)$ from
above is given by the separable symmetric Frobenius algebra
$A(Q) = \bigoplus_{g\in G} I^g_V$, from which one can extract complete
information about the orbifolded topological Landau-Ginzburg model. 

\vspace{0.1cm}
The most interesting orbifold equivalences so far have been found for simple
singularities of ADE type \cite{CR3,CRCR}. The potentials are
$$
V_{A_n}= x_1^{n+1}+x_2^2\ ,\ \ 
V_{D_{d}}=x_1^{d-1}+x_1x_2^2\ ,\ \ 
V_{E_6}=x_1^3+x_2^4\ ,\ \
V_{E_7}=x_1^3+x_1x_2^3\ ,\ \ 
V_{E_8}=x_1^3+x_2^5
$$
with $d\geq 2$ and $n\geq 4$, the corresponding Landau-Ginzburg models
are related to so-called $N=2$ superconformal minimal models with central
charge $\hat c <1$. It turns out that whenever two of these
potentials have the same central charge, they are also orbifold equivalent; the
classes with more than one representative are $\{A_{d-1},D_{d/2+1}\}$ for even
$d$ not equal to 12, 18 or 30, and $\{A_{11}, D_7 ,E_6\}$,
$\{A_{17}, D_{10}, E_7\}$ and $\{A_{29}, D_{16}, E_8\}$. The $A$-$D$ orbifold
equivalences are related to (simple current) orbifolds in the CFT context,
but the $A$-$E$ orbifold equivalences do not arise from any group action
\cite{CRCR}; they are examples of ``symmetries'' beyond groups.

For the purposes of elucidating some general observations to be made later,
and also of conveying an idea of the typical complexity of the matrix
factorisations
involved, we reproduce a concrete example of an orbifold equivalence from
\cite{CRCR}, namely that between $V_{A_{11}} = x_1^{12} + x_2^2$ and
$V_{E_6} = y_1^3 + y_2^4$. In this case, the smallest possible
(see subsection 3.2) orbifold equivalence 
is of rank 2. With $Q$ formed from $E$ and $J$ as in (\ref{QEJ}),  
the matrix elements of $E$ are given by
\bea \la{EforA11E6}
E_{11}  &=&\vphantom{\sum}  y_2^2-x_2 +  {\textstyle \frac{1}{2}} y_1(sx_1)^2 +
{\textstyle\frac{2t+1}{8}} (sx_1)^6     
  \cr
  E_{12}  &=&\vphantom{\sum}  -y_1 + y_2(sx_1) +    {\textstyle \frac{t+1}{4}} (sx_1)^4           
\cr
E_{21}  &=&\vphantom{\sum}  y_1^2 + y_1y_2(sx_1) +    {\textstyle \frac{t}{4}} y_1(sx_1)^4
 +  {\textstyle \frac{2t+1}{4}}y_2 (sx_1)^5 -  {\textstyle \frac{9t+5}{48}}(sx_1)^8    
 \cr
 E_{22}  &=&\vphantom{\sum}  y_2^2+x_2 +   {\textstyle \frac{1}{2}} y_1(sx_1)^2
   +  {\textstyle \frac{2t+1}{8}}  (sx_1)^6      
\eea
$J = - \mathrm{adjugate}(E)$, and the complex coefficients $s,t$ satisfy the
algebraic equations \hbox{$t^2 = 1/3\ ,\ \  s^{12} = -576(26t-15)$}.
This defect has non-zero quantum dimensions, namely $q_L(Q) = s,\ q_R(Q) = 3(1-t)/s$.

In \cite{NRC1}, an explicit rank 4 orbifold equivalence was written down
between two of the fourteen (quasi-homogeneous) exceptional unimodal
Arnold singularities (list e.g.\ in \cite{ArnoldBook}), namely between
$E_{14}$ and $Q_{10}$ described by the potentials
$$
V_{E_{14}}(x) = x_1^8 + x_2^2 + x_3^3 
\ \ \mathrm{and}\ \
V_{Q_{10}}(y) = y_1^4 + y_1 y_2^2 + y_3^3\ , 
$$
both having central charge $\hat c = \frac{13}{12}$. However, one
should notice that this orbifold equivalence already follows from
the $A$-$D$ results of \cite{CR3} and the general property Theorem 2.2 (c):
One can write $E_{14} = A_7 \times A_2$ and $Q_{10} = D_5 \times A_2$, 
and one has $A_7 \oeq D_5$. In the same way, one can of course construct
other orbifold equivalences at arbitrarily high central charge,
simply by ``adding up'' suitable simple singularities with $\hat c <1$.

\section{Some structural results on orbifold equivalences}\label{sec:GP}

If one tries to generate examples of orbifold equivalences truly beyond
simple singularities, one soon realises that the approach taken in \cite{CRCR}
is neither general nor systematic enough. In that work, the method employed
to find expressions like (\ref{EforA11E6}) was to set one of the variables
$x_i, y_j$ occurring in $W(x,y) = V_1(x)-V_2(y)$ to zero, to pick some
simple matrix factorisation $\widetilde{Q}$ of the resulting potential
$\tilde W$ and to complete $\widetilde{Q}$ to a graded matrix factorisation
$Q(x,y)$ of the full $W(x,y)$ using quasi-homogeneous entries that contain
the missing variable -- under additional simplifying constraints such
as $J = - \mathrm{adjugate}(E)$. But as soon as one has to cope with a
larger number of variables, or higher rank, one needs a lot of luck 
to hit a good starting point $\widetilde{Q}$.

Nevertheless, the computations in \cite{CRCR} contain germs of ideas which 
can be formulated in general terms and exploited in a systematic manner.
In this section, we will show that 
every graded orbifold equivalence has a (finite) perturbation expansion,
a structure from which one can draw some general conclusion on the form $Q(x,y)$
must take. The grading is a crucial ingredient, and we will present
a constraint on the allowed grading matrices in subsection 3.2. Together
with the perturbative structure, this will enable us to devise a relatively
efficient search algorithm for orbifold equivalences in section 4.

\subsection{Orbifold equivalences as graded perturbations}

Given a matrix factorisation $Q$ of a potential $W$, one can ask whether
$Q$, or the differential $d_{QQ}$ associated with it, admits deformations.
As is familiar in the context of deformation theory, deformation
directions are controlled by $\mathrm{Ext}^1$ -- or $H^1_Q$, the space
of boundary fermions --, obstructions by $\mathrm{Ext}^2$ -- or $H^0_Q$.
References and some results can e.g.\ be found in \cite{CDR}.

We will now show that graded orbifold equivalences, or indeed any
graded defect between $V_1(x)$ and $V_2(y)$, can be naturally viewed
as a deformation of a matrix factorisation of $V_1(x)$, with the
variables $y_j$ featuring as deformation parameters and
$-V_2(x) \mathbf{1}_{2N}$ as obstruction term. Since we are not
restricting to the vanishing locus of the obstruction term, we
speak of ``perturbations'' rather than deformations.

That the defect is graded has a very desirable consequence: the
perturbation expansion terminates after finitely many steps.

\vspace{.3cm}
For this discussion, it is convenient to introduce some further 
notions concerning graded matrix factorisations:
Let $U(\lambda)$ be a grading matrix as in (\ref{grading}). 
Borrowing some further physics terminology, we say that 
a matrix $A\in M_{2N}(\C[z])$ has ``R-charge'' $R$ wrt.\ the
grading $U(\lambda)$ if
\be\la{rchargedef}
U(\lambda)\, A(\lambda \rhd z)\, U(\lambda)^{-1}  = \lambda^R\,A(z)\ .
\ee
(The graded matrix factorisation $Q$ in (\ref{grading}) then has R-charge 1.) 

This relation implies that the entries $A_{rs}$ of $A$ are quasi-homogeneous
polynomials in the $z_i$. Their weights can be computed from the grading
matrix as
\be
w(A_{rs}) =  g_s-g_r+R\quad\mathrm{for}\ \ r,s = 1,\ldots,2N\ . 
\ee
In the special case $A=Q$ we will sometimes write $w(Q)$ for the matrix 
formed from the $w(Q_{rs})$ and call it the weight matrix of $Q$; analogously
we will use $w(E)$, $w(J)$ for the weight matrices of
$E$ and $J$ related to $Q$ as in (\ref{QEJ}).

\vspace{.2cm}
By way of a brief excursion, and also as a step towards a proof of
Prop.\ 2.3, let us use the notion of R-charges to provide a self-contained
derivation of a statement that is well-known in the physics literature
on topological Landau Ginzburg models, namely that the correlators
(\ref{bulkcorr},\ref{bdycorr}) have a ``background charge''. 
Instead of employing arguments from an underlying twisted conformal field
theory, this can be derived from properties of the residue. We focus on
the Kapustin-Li correlator here: 

\vspace{.2cm}\noindent
\textbf{Proposition 3.1: } Set $z = (z_1,\ldots,z_k)$ and let $Q(z)$ be a 
(graded) rank $N$ matrix factorisation of a potential $W(z)$ with
dim$_{\C}(\mathrm{Jac}(W)) < \infty$.
Let $\psi \in M_{2N}(\C[z])$ be a morphism of definite $\Z_2$-degree
$s(\psi)$ and definite R-charge $R_\psi$. Then 
$$
\langle\, \psi\,\rangle_{Q}^\mathrm{KapLi} = 0
$$
unless $s(\psi) + k$ is even and unless $R_\psi = \hat c(W)$. 

\vspace{.2cm}\noindent
\textsl{Proof:}
The statement on the $\Z_2$-degree follows because $Q$ and its partial
derivatives are odd matrices wrt.\ to the $\Z_2$-grading $\sigma$,
hence a product of $n$ of these with the even or odd matrix $\psi$
has no diagonal terms, hence zero supertrace, if $k+s(\psi)$ is odd.

As the Jacobi ring of $W$ is a finite-dimensional $\C$-vector space,
for each $i=1,\ldots,k$ there is a $\nu_i\in \Z_+$ and polynomials $C_{ij}$
such that $z_i^{\nu_i} = \sum_j C_{ij}(z) \partial_{z_j}W(z)$. This implies,
see e.g.\ \cite{GriffithsHarris}, that 
$$
  \mathrm{res}_z \biggl \lbrack
       { f  \over \partial_{z_1}W \cdots \partial_{z_k}W } \biggr \rbrack = 
     \mathrm{res}_z \biggl \lbrack
       { \det(C) f  \over z_1^{\nu_1} \cdots z_k^{\nu_k} }\biggr \rbrack
$$
for any polynomial $f(z)$.
       
In the case at hand,
$f = \mathrm{str}(\partial_{z_1}Q\cdots\partial_{z_k}Q\cdot\psi)$, and
since $Q$, its derivatives, and $\psi$ have definite R-charges, 
a rescaling of the $z_i$ can be traded for conjugation
with the grading matrix $U(\lambda)$
-- this leaves the supertrace invariant -- up to extra prefactors
$\lambda^{1-|z_i|}$ resp.\ $\lambda^{r\psi}$ from relation (\ref{rchargedef}).
Hence $f$ is quasi-homogeneous of weight $R_\psi + \hat c(W)$. 

It is easy to see that $\det(C)$ is quasi-homogeneous of weight
$-k - \hat c(W) + \sum_i\,\nu_i|z_i|$. The residue projects $f\cdot \det(C)$
onto the monomial $z_1^{\nu_1-1} \cdots z_k^{\nu_k-1}$, which has weight
$-k + \hat c(W) + \sum_i\,\nu_i|z_i|$. Thus the residue can be non-zero
only if $R_\psi = \hat c(W)$.    $\quad\square$

\vspace{.3cm}
Turning to the perturbation expansion of orbifold equivalences, we assume,
as before, that $V_1(x)$ and $V_2(y)$ are quasi-homogeneous potentials of
weight 2, without linear terms, and we denote the weights of the variables by
$|x_i|$ for $i=1,\ldots,n$ resp.\ $|y_j|$ for $j=1,\ldots,m$.
We abbreviate $W(x,y) := V_1(x)-V_2(y)$. 

\vspace{.2cm}\noindent
\textbf{Proposition 3.2: } Assume that $Q(x,y)$ is a (graded) rank $N$ orbifold
equivalence between $V_1(x)$ and $V_2(y)$, i.e.\ $Q^2 = W\,\mathbf{1}_{2N}$ 
and $q_L(Q)\,q_R(Q) \neq  0$.
Set $Q_1(x) := Q(x,y)|_{y=0}$ and
$F_j := \partial_{y_j}Q(x,y)|_{y=0}$ for $j=1,\ldots,m$. Then

\begin{itemize} 

\item[(1)]\ \ $F_j$ is a fermionic morphism of $Q_1$ with R-charge $R_j = 1-|y_j|$, for all
$j=1,\ldots,m$.

\item[(2)]\ \ The left quantum dimension of $Q$ can be written as 
 $$q_L(Q) = \langle\, F_1 \cdots F_m\,\rangle_{Q_1}^\mathrm{KapLi}
 $$
where $\langle\cdots\rangle_{Q_1}^\mathrm{KapLi}$ denotes the Kapustin-Li boundary correlator
of the LG model with bulk potential $V_1$ and boundary condition $Q_1$.

\item[(3)]\ \ $Q(x,y)$ has a finite perturbation expansion, with $y_j$ appearing as parameters:
 $$ 
 Q(x,y) = \sum_{\kappa=0}^{\kappa_\mathrm{max}} Q^{(\kappa)}(x,y)  
 \quad \mathrm{ with}\ \ Q^{(0)}(x,y)= Q_1(x)\ \mathrm{and}\ Q^{(1)}(x,y) = \sum_j y_j\,F_j\,.
 $$
The higher
 order terms satisfy
 \be \la{higherpert}
 \{ Q_1,Q^{(\kappa)}\} + \sum_{\lambda=1}^{\kappa-1} Q^{(\kappa-\lambda)} Q^{(\lambda)}
      = -V_2^{(\kappa)}
 \ee
 where $\{\cdot,\cdot\}$ denotes the anti-commutator and where
 $V_2^{(\kappa)}$ is the order $\kappa$ term of $V_2$. 
 
\end{itemize}

\vspace{.2cm}\noindent
\textsl{Proof:} The $F_j$ are odd wrt.\ the $\sigma$-grading (as $Q$ is), and they are
in the kernel of $d_{Q_1Q_1}$ because $\{ Q_1, F_j\} = \partial_{y_j}Q(x,y)^2|_{y=0} =
\partial_{y_j}W\,\mathbf{1}|_{y=0} = 0$. Let $U(\lambda)$ be the grading for $Q(x,y)$ --
and, for that matter, for $Q_1(x)$. Differentiating (\ref{grading}) gives 
$$
U(\lambda)\, \partial_{y_j}Q(\lambda \rhd x, \lambda \rhd y)\, U(\lambda)^{-1}  =
\lambda^{1-|y_j|}\,\partial_{y_j}Q(x,y)\ ,
$$  
so in particular the $F_j$ are fermions of $Q_1$ with definite R-charge
$R_j = 1 - |y_j|$.

To see the second statement, note that the left quantum dimension of a
graded matrix factorisation is a quasi-homogeneous polynomial in $y$, 
and in fact has to be a (non-zero) number in order to be invertible.
Hence $q_L(Q)$ does not depend on the $y$-variables. 
Setting $y=0$ in the first of the residue formulas (\ref{qdimdef}) directly
produces the Kapustin-Li correlator  (\ref{bdycorr}) of the
product $F_1\cdots F_m$ of ``boundary fermions'' in the $(V_1,Q_1)$ theory.

Eq.\ (\ref{higherpert}) simply follows from a Taylor expansion of
$Q^2 = W\,\mathbf{1}_{2N}$ around $y=0$, keeping in mind that $Q_1(x)$ is a
matrix factorisation of $V_1(x)$.

Finiteness of the perturbation series can be seen by analysing R-charges and
weight matrices: The order $\kappa$ term
$Q^{(\kappa)}(x,y) = \sum_{\vec p} M^{(\kappa)}_{\vec p} y^{\vec p}$ is a linear combination
of monomials in the $y_j$, $\vec p \in \Z_+^{m}$ with $p_1+\ldots+p_m = \kappa$,
with matrix-valued coefficients $M^{(\kappa)}_{\vec p} \in M_{2N}(\C[x])$. These
coefficients are odd wrt.\ the $\mathbb{Z}_2$-grading, and they have
R-charge $R^{(\kappa)}_{\vec p} = 1 - p_1|y_1| - \ldots - p_m|y_m|$.
The entries of $M^{(\kappa)}_{\vec p}$ are homogeneous polynomials in the $x_i$, the 
weight of the $r$-$s$-entry is $w_{rs} = g_s -g_r+ R^{(\kappa)}_{\vec p}$, where the
$g_r$ define the grading matrix $U(\lambda)$ as before.
Since all the variable weights are strictly positive, for large enough $\kappa$
not only the R-charge but also the weights $w_{rs}$ will become negative for
all $r,s \in \{1,\ldots,2N\}$, which implies that $M^{(\kappa)}_{\vec p}$ has to vanish. 
$\quad\square$ 

\vspace{.2cm}

An analogous expansion can be performed around $x=0$, and the right quantum
dimension of $Q(x,y)$ takes the form of a correlator of boundary fermions
in the Landau-Ginzburg with bulk potential $-V_2(y)$ and boundary condition
$Q_2(y)) := Q(x,y)|_{x=0}$. Setting 
$\tilde F_i := \partial_{x_i}Q(x,y)|_{x=0}$, we have
$$
 q_R(Q) = \langle\, \tilde F_1 \cdots \tilde F_n\,\rangle_{Q_2}^\mathrm{KapLi} \ .
$$

\vspace{.1cm}
In the present paper, the main application of Prop.\ 3.2 will be to devise
a more systematic search algorithm for orbifold equivalences, which allows us
to tackle more difficult situations than the simple singularities discussed
in \cite{CRCR}. But there are some immediate structural consequences implied by
the perturbation expansion:

\vspace{.1cm}

First off, the condition $\hat c(V_1) = \hat c(V_2)$ necessary for the existence
of a graded orbifold equivalence follows immediately from the ``background
charge'' of topological Landau-Ginzburg correlators: The product
$F_1\cdots F_m$ is a morphism
with R-charge $\sum_j(1-|y_j|) = \hat c(V_2)$, and its Kapustin-Li correlator
in the $(V_1,Q_1)$-model vanishes, according to the statement rederived in
Prop.\ 3.1, unless this R-charge coincides with the background charge $\hat c(V_1))$
of that model.

\vspace{.2cm}
Prop.\ 3.2 also constrains what form the matrix elements of an orbifold
equivalence $Q(x,y)$ can take: Clearly, for each $j=1,\ldots,m$ there
must be a $Q$-entry that contains a term linear in $y_j$,  
lest one of the partial derivatives $F_j$ is zero; likewise for the $x_i$.
(In fact, none of the $F_j$ can be trivial in the $Q_1$-cohomology, i.e.\ none
can be of the form $F_j = [Q_1,A_j]$ for some $A_j$, because the Kapustin-Li
form is independent of the representative of the cohomology class.)

\vspace{.1cm}
Moreover, under very mild additional assumptions on the potentials, one can
show that orbifold equivalences must involve some ``entanglement'' of the
$x$- and $y$-variables:

\vspace{.2cm}\noindent
\textbf{Proposition 3.3: } Assume that $V_2(y) \in \maxm^3$, i.e.\ has no quadratic
or lower order terms. Then an orbifold equivalence $Q(x,y)$ between $V_1(x)$
and $V_2(y)$ must have mixed $xy$-terms, i.e.\ it cannot have the form
$Q(x,y) = Q_1(x) + Q_2(y)$. 

\vspace{.1cm}\noindent
\textsl{Proof: } Assume $Q(x,y) = Q_1(x) + Q_2(y)$. The first summand is a matrix
factorisation of $V_1(x)$, the second summand one of $-V_2(y)$; consequently
$\{Q_1(x),Q_2(y)\} = 0$ and also $\{\partial_{x_i}Q_1, F_j\} = 0$. Moreover,
$$
 0=  - \partial_{y_{j_1}} \partial_{y_{j_2}} V_2(y)|_{y=0} =
  \{ F_{j_1}, F_{j_2}\}
  + \{ Q_2(y), \partial_{y_{j_1}} \partial_{y_{j_2}}Q_2(y)\}|_{y=0}\ .
$$ 
The last term vanishes (since the matrix factorisations we consider have
no constant terms), so all the $F_j$ anti-commute and square to zero.

Let $N := \sigma\,\partial_{x_1}Q_1\cdots \partial_{x_n}Q_1\,F_1\cdots F_m$,
which is the argument of the trace in the residue formula for $q_L(Q)$.
In the case at hand, this matrix $N$ is nilpotent, 
$$
  N^2 = \pm (\partial_{x_1}Q_1\cdots \partial_{x_n}Q_1)^2(F_1\cdots F_m)^2 = 0\ \ ,
$$
hence $\mathrm{tr}(N)=0$ and $q_L(Q) = 0$.

Note that we can relax the assumption on $V_2(y)$: as soon as there is one
variable $y_{j^*}$ such that $\partial_{y_{j^*}}\partial_{y_j}V_2(y)|_{y=0} = 0$
for all $j=1,\ldots,m$, we have that $F_{j^*}$ anti-commutes with all
$F_j$, and $N$ is nilpotent.    $\quad\square$

\vspace{.2cm}
 
In particular, this result rules out the simplest tensor products as orbifold
equivalences (under the stated assumptions on the potentials): 
if $Q(x,y) = Q_a(x)\, \hat\otimes\, Q_b(y)$ where
$Q_a(x)$ is a matrix factorisation of $V_1(x)$ and $Q_b(y)$ one of $-V_2(y)$,
then $Q(x,y)$ has zero quantum dimensions.

\vspace{.1cm}
That the standard method (forming tensor products) of constructing matrix
factorisations for complicated polynomials is barred when seeking orbifold
equivalences goes some way in explaining why the latter are hard to find.
Results of the type of Prop.\ 3.3 may also prove useful for showing that
equality of central charges is an insufficient criterion for
two potentials to be orbifold equivalent.

\subsection{Weight split criterion}

The perturbative expansion described in Prop.\ 3.2 is a useful ingredient 
of an algorithmic search for orbifold equivalences, but, as it stands,
the need to select a grading and a $Q_1(x)$ as starting point seems to 
limit efficiency quite severely. In this subjection, we will
point out that the gradings (i.e.\ the weight matrices of $Q_1(x,y)$
and $Q(x,y)$) are subject to a highly selective criterion -- a criterion
that applies to any graded matrix factorisation $Q(z)$ of any
quasi-homogeneous potential $W(z)$, not just to defects. 

It will be more convenient to rescale the variable weights such
that all $|z_i|$ are natural numbers; so for the time being, the
weight of $W(z)$ is given by some integer $D_W\in\Z_+$, not necessarily
equal to 2.

\vspace{.2cm} 
Before giving a general formulation, let us see the criterion ``at work''
in the concrete example of the $A_{11}$-$E_6$ orbifold equivalence found
in \cite{CRCR} and reproduced in subsection 2.2.
Here, $V_{A_{11}}(x) = x_1^{12} + x_2^2$, $V_{E_6}(y) = y_1^3 + y_2^4$,
and $W(z) = V_{A_{11}}(x) - V_{E_6}(y)$ with $z=(x,y)$.
The variable weights are $|x_1 | = 1, |x_2 | = 6, |y_1 | = 4, |y_2 | = 3$
(after scaling up to integers, so that $D_W=12$.)

Any graded matrix factorisation $EJ = JE = W \mathbf{1}_{N}$  must in particular
contain (quasi-homogeneous) polynomials factorising the $x_2^2$-term from $W$
-- and such factors must occur in each row and each column of $E$ and $J$.
Up to constant prefactors, these polynomials must be of the form 
$x_2 + f_{rs}$ for some $f_{rs}$ having the same weight as $x_2$. 
So each row and each column of the weight matrices $w(E)$ and $w(J)$ must
contain a $6$.

Likewise, the $y_1^3$-term has to be factorised, so each row and column
of $w(E)$ and $w(J)$ has to contain a $4$ (from a factor $y_1^1+\ldots$) or an
$8$ (from a factor $y_1^2+\ldots$).
 
If we want to construct a rank $N=2$ matrix factorisation of $W = 
V_{A_{11}} - V_{E_6}$, these two observations (together with the constraint
that $Q$ should be graded) fix the weight matrices  completely, up to
row and column permutations and up to swapping $E$ and $J$: 
  \bea
  w(E)&=& 
   \left(
  \begin{array}{cc}
			6 & 4 \\
 			8 & 6
  \end{array}
  \right)
  \eea
which is indeed the weight matrix for the $A_{11}$-$E_6$ orbifold equivalence
(\ref{EforA11E6}) found by Carqueville et al. (Thanks to the low rank and
the small number of variables, it is fairly easy to arrive at a concrete
$Q$ once the above $w(E)$ is known.) 

\vspace{.2cm}
In order to formulate the criterion in general, we need some notation. Let
$$
W(z) = \sum_{\tau=1}^T m_\nu(z)
$$
be the decomposition of the potential into monomial terms;
each $m_\tau$ has weight $D_W$.
For each $\tau=1,\ldots,T$, let $S_\tau$ be the set of weights of possible
non-trivial divisors of $m_\tau$, i.e.\
$$
S_\tau = \Big\{\; w \in \{1,\ldots,D_W-1\}\; :\; \exists f \in \C[z]\ \;
\mathrm{s.th.}\ 
f\ \mathrm{divides}\ m_\tau\ \mathrm{and}\ f\  \mathrm{has\ weight}\ w\; \Big\}
$$

\vspace{.1cm}\noindent
\textbf{Weight split criterion: } If $Q(z)$ is a graded matrix factorisation of $W(z)$ 
with weight matrix $w(Q)$, then each row and and column of $w(Q)$ contains an
element of $S_\tau$ for all $\tau=1,\ldots,T$.

\vspace{.2cm}
Let us look at two further examples to illustrate the usefulness of this criterion.
For the two unimodal Arnold singularities 
$V_{E_{13}}(x)  = x_2^3 + x_2x_1^5$ and $V_{Z_{11}}(y) =  y_1^3y_2 + y_2^5$, 
the variable weights are $|x_1 | = 2, |x_2 | = 5, |y_1 | = 4, |y_2 | = 3$
(re-scaled so that $W$ has weight $15$). 
The terms in $W$ admit weight splits $5+10$ and  $5+10=7+8=9+6=11+4=13+2$
(from $E_{13}$)  and $12+3=8+7=4+11$ and $3+12=6+9$  (from $Z_{11}$). In each
row and each column of $w(E)$, there must be a $5$ or a $10$, and
there must be one from the set $\{3,4,7,8,11,12\}$. 

One can just about fit the above weights into a rank 2 matrix
$w(E)$ with entries $5, 12, 10, 3$, but this leads to zero quantum
dimensions (the associated $Q$ are tensor products and ruled out as
orbifold equivalences by Prop.\ 3.3).

At rank 3, one can form 24 weight matrices $w(E)$ satisfying the weight split
criterion, and one of those leads to an orbifold equivalence, see the next
section. It is worth mentioning that the ``successful'' $w(E)$ is one where
many entries are members of \textsl{both} the weight split lists $S_\tau$ coming from
$E_{13}$ \textsl{and} the weight split lists coming from $(Z_{11})$; these
offer the best opportunity for an ``entanglement'' of $x$ and $y$ variables. 

How restrictive the weight split criterion can be is demonstrated when one tries
to construct an orbifold equivalence for the Arnold singularities $Z_{13}$ and
$Q_{11}$, see section 4: here, one needs a rank 6 matrix factorisation, and of
about 2.7 million
conceivable weight matrices $w(Q)$ only 60 pass the criterion.

\vspace{.2cm}
There are additional restrictions on viable weight matrices $w(Q)$ which
apply if $Q$ is to be an orbifold equivalence between $V_1(x)$ and
$V_2(y)$. E.g., the requirement
that non-trivial fermions of given R-charge have to exist (needed for
non-zero quantum dimensions, cf.\ Prop  3.2), means that for each
variable $y_j$, at least one of the $w(Q)$-entries must be of the
form $|y_j|+n(x)$ where $n(x)$ is some $\Z_+$-linear combination of the
weights $|x_i|$; analogously with the roles of $x$ and $y$ interchanged.
In the examples we studied, this condition from existence of fermions turned
out to be far less restrictive than the weight split criterion arising
from the matrix factorisation conditions.

\section{Algorithmic search, and some concrete results}\label{sec:examples}

In this chapter, we will present some new examples of orbifold equivalences.
Most of them were discovered using an algorithm based on the perturbative expansion
introduced in the previous section. First, we make some general remarks on
the ``computability'' of orbifold equivalences and outline a
computer-implementable algorithm
to deal with the problem, then we list the new examples themselves. 

\subsection{Towards an algorithmic search for orbifold equivalences}

The question whether there is a rank $N$ orbifold equivalence $Q$ between two
given potentials $V_1$ and $V_2$ can be converted into an ideal membership
problem and, for fixed $N$, can be decided by a finite computation.

To see this, let us write the matrix elements of $Q$ as
\be\la{qmels}
Q_{rs} = \sum_{\vec p}\; a_{rs,\vec p}\,z^{\vec p}
  \quad \mathrm{for}\ \ r,s \in \{1,\ldots,2N\}
\ee
where $z=(x_1,\ldots,x_n,y_1,\ldots,y_m)$ and where $\vec p \in \Z_+^{m+n}$ is a
multi-index. The main ``trick'' now is to shift one's focus away from the
variables $z$ and work in a ring of polynomials in the $a_{rs,\vec p}$: 

The requirement that $Q$ is a rank $N$ matrix factorisation of
$W(z) = V_1(x)-V_2(y)$ imposes polynomial (in fact:\ bilinear) equations
$f_\alpha^\mathrm{MF}(a) = 0$ on the coefficients $a_{rs,\vec p} \in \C$.
($\alpha$ labels the various bilinear equations, $a$ collectively denotes all
the coefficients.)

The quantum dimensions can be computed, using definition (\ref{qdimdef}),
whether or not $Q$ is a matrix factorisation; for a graded $Q$, one obtains
two polynomials (of degree $n+m$) in the $a_{rs,\vec p}$. The requirement
that both quantum dimensions are non-zero is equivalent to the single
equation
$$
f^\mathrm{qd}(a) :=  q_L(Q) q_R(Q)\, a_\mathrm{aux} -1 = 0
$$
being solvable, where $a_\mathrm{aux}$ is an additional auxiliary coefficient.

Thus, the matrix $Q$ is an orbifold equivalence between $V_1$ and $V_2$
if and only if the system 
\be \la{oeqeqs}
f_\alpha^\mathrm{MF} = 0\ , \quad f^\mathrm{qd} = 0
\ee
of polynomial equations in the coefficients  $a_{rs,\vec p}$ and $a_\mathrm{aux}$ has
a solution. By Hilbert's weak Nullstellensatz, this is the case iff
\be \la{idealmemb} 
1 \notin \langle\, f^\mathrm{MF}_\alpha, f^\mathrm{qd}\,\rangle_{\C[a]} \ .
\ee
This type of ideal membership problem can be tackled rather efficiently
with computer algebra systems like Singular. (Such systems are usually
restricted to working over $\Q$, but for potentials $V_1,\;V_2$ with
rational coefficients it is enough to study (\ref{idealmemb}) over the
rationals in order to prove or disprove existence of an orbifold equivalence
with coefficients $a_{rs,\vec p}$ in the algebraic closure $\overline{\Q}$.) 

\vspace{.2cm}
Once a grading $U(\lambda)$, hence a weight matrix for $Q$, has been chosen,
it is easy to write down the most general homogeneous matrix elements
$Q_{rs}$ (\ref{qmels}) that conform with this grading. Moreover, there is
only a finite number of possible gradings
$U(\lambda) = \mathrm{diag}(\lambda^{g_1},\ldots,\lambda^{g_{2N}})$ for a given
rank $N$. To see this, recall that the weights of the $Q$-entries are given
by $w(Q_{rs}) = g_s - g_r +1$ (we set the weight of the potential to 2),  
and also that we can fix $g_1=0$ wlog -- so in particular $w(Q_{1r}) = g_r +1$ 
and $w(Q_{r1}) =  - g_r +1$. Therefore, at least one of the $g_r$ has to
satisfy $-1\leq g_r \leq 1$, otherwise the entire first row
or column of $Q$ would have to vanish (because the weights would all be negative),
which would contradict the matrix factorisation conditions. We can repeat
the argument for the $g_r$ nearest to $g_1$ and find, overall, that
$g_r \in \lbrack -2N,2N\rbrack$ for all $r$. Finally, $Q_{rs}$ can be a
non-zero polynomial in the $x_i,\,y_j$ only if its weight $w(Q_{rs})$ is a
sum of the (finitely many, rational) weights $|x_i|,\,|y_j|$, hence only
finitely many choices $g_r$ from the interval $\lbrack -2N,2N\rbrack$ can
lead to a graded rank $N$ matrix factorisation of $V_1(x)-V_2(y)$.

All in all, the question whether there exists a rank $N$ orbifold equivalence
between two given potentials $V_1,\,V_2$ can be settled in principle.
Our guess is that there is an upper bound $N_\mathrm{max}(V_1,V_2)$ such that,
if no orbifold equivalence of rank $N < N_\mathrm{max}(V_1,V_2)$ exists, then
none exists at all -- but we have only circumstantial evidence: all
known (indecomposable) examples of orbifold equivalences have rank smaller
than the nested tensor product matrix factorisation obtained by
factorising each monomial
in $V_1-V_2$; and packing a matrix factorisation ``too loosely'' risks making
the supertrace inside the quantum dimensions vanish.

\vspace{0.3cm}
So much for the abstract question whether orbifold equivalence is a property
that can be decided algorithmically at all. In order to search for concrete
examples, we have devised an \textsl{algorithm} based on the perturbation expansion
and the weight split criterion introduced in section 3:

\begin{itemize} 
  
\item[(a)] From the potentials $V_1(x),\,V_2(y)$, compute the variable weights
$|x_i|,\,|y_j|$.

\item[(b)] Choose a rank $N$.

\item[(c)] Exploiting the weight split criterion from subsection 3.2, compute all
admissible gradings (i.e.\ weight matrices) for this rank. 

\item[(d)] Choose a weight matrix and form the most general matrix factorisation
$Q_1(x)$ of $V_1(x)$ with this weight matrix. 

\item[(e)] For each $y_j$, compute the space of fermions $F_j$ of $Q_1(x)$ with
R-charges $1-|y_j|$. 

\item[(f)] For any R-charge $R_M$ that can occur in the expansion of $Q(x,y)$ from
Prop.\ 3.2, determine the space of odd matrices with that R-charge.

\item[(g)] Compute $Q(x,y)$ using the conditions from Prop.\ 3.2 (c), then
  compute the quantum dimensions $q_L(Q)$ and $q_R(Q)$. (Everything will depend on
  unknown coefficients $a$.) 

\item[(h)] Extract the conditions $f^\mathrm{MF}_\alpha(a)=0$ and $f^\mathrm{qd}(a)=0$ on the
coefficients appearing in $Q(x,y)$ and check whether this system of
polynomial equations admits a solution. 

 \end{itemize}

\vspace{0.2cm}

Computer algebra systems such as Singular have in-built routines to perform
the last step, employing (variants of) Buchberger's algorithm to compute
a Gr\"obner basis of the ideal spanned by $f^\mathrm{MF},\,f^\mathrm{qd}$.

Already when forming $Q_1(x)$ with a given weight matrix, undetermined
coefficients $a$ enter the game -- but far fewer than would show up in
the most general matrix $Q(x,y)$ with the same weight matrix, because
one only uses the $x$-variables to form quasi-homogeneous entries:  
the perturbation expansion from Prop.\ 3.2 ``organises'' the computation
to some extent from the outset. Nevertheless, even for harmless looking
potentials $V_1,\;V_2$ one can easily end up with
close to one thousand polynomial equations in hundreds of unknowns
$a_{rs,\vec p}$. Due to restrictions on memory and run-time, it is advisable
in practice to make guesses for some of the coefficients $a_{rs,\vec p}$
occurring in $Q(x,y)$ or already in $Q_1(x)$, instead of trying to tackle
the most general ansatz. We have succeeded in automatising most of
the steps involved in making the equations tractable for Singular, some
of the results are collected in the next subsection.

Finding an explicit solution for the coefficients $a$ is of course desirable,
but not necessary to prove that two potentials are orbifold equivalent. 
It appears that Singular is not the optimal package for determining
explicit solutions (although it is very efficient in establishing solvability);
feeding the polynomial equations resulting from the Singular code into
Mathematica, say, might be more promising.

\vspace{0.1cm}
If one is content with existence statements, additional avenues are open:
One could e.g.\ employ numerical methods to find approximate solutions to the
system of equations (\ref{oeqeqs}), then check whether any of them satisfies
the criteria of the Kantorovich theorem or of Smale's $\alpha$-theory. If so,
one has proven (rigorously) that there is an exact solution in a
neighbourhood of the numerical one.
We did not take this route, but it might lead to a more efficient
computational tool towards a classification
of orbifold equivalent potentials.

\subsection{New examples} 

We now present new examples of orbifold equivalences, starting with a few
isolated (but hard-won) cases, including all remaining pairs of
unimodal Arnold singularities, then adding a series of equivalences obtained
by simple transformations of variables.

\vspace{0.2cm}\noindent
\textbf{Theorem 4.1: } In each of the following cases, the potential
$V_1(x)$ is orbifold equivalent to the potential $V_2(y)$:

\begin{itemize}

\item[(1)] $V_1(x) = x_1^6 + x_2^2$ and $V_2(y) = y_1^3 + y_2^3$.

(These are the singularities $A_5$ resp.\ $A_2 \times A_2$, at central charge
$\hat c = {2\over3}$.)  

\item[(2)] $V_1(x) = x_1^5 \,x_2 + x_2^3$ and $V_2(y) = y_1^3 \,y_2 + y_2^5$. 

(These are two of the exceptional unimodal Arnold singularities,
namely $E_{13}$ resp.\ $Z_{11}$, at central charge $\hat c = {16\over15}$.)

\item[(3)] $V_1(x) = x_1^6+x_1 \,x_2^3+x_3^2$  and 
$V_2(y) =  y_2 \,y_3^3 + y_2^3 + y_1^2 \,y_3$. 

(These are the exceptional unimodal Arnold singularities 
 $Z_{13}$ resp.\ $Q_{11}$, at central charge $\hat c = {10\over9}$.)

\item[(4)] $V_1(x) = x_1^2 \,x_3 + x_2 \,x_3^2 + x_2^4$  and 
$V_2(y) =  -y_1^2+ y_2^4 + y_2 \,y_3^4$. 

(These are the exceptional unimodal Arnold singularities 
 $S_{11}$ resp.\ $W_{13}$, at central charge $\hat c = {9\over8}$.)

\item[(5)] $V_1(x) = x_1^{10} x_2+x_2^3$ and $V_2(y) = y_1 y_2^7+ y_1^3 y_2$.
  
  (These are a chain resp.\ a loop (or cycle), in the nomenclature of
  \cite{KrSk,HertlingKurbel},  at central charge
  $\hat c = {6\over5}$, a value shared by the pair $Q_{17}$ and $W_{17}$ 
  of bimodal Arnold singularities.) 
  
\end{itemize}

\vspace{.2cm} 
\noindent
\textsl{Proof:} In contrast to $E_{14}$-$Q_{10}$, none of these 
cases can be traced back to known results on simple singularities. 
Lacking, therefore, any elegant abstract arguments, we can only establish
these orbifold equivalences by finding explicit matrix factorisations
$Q$ of $V_1-V_2$ with non-zero quantum dimensions. The ranks of the $Q$
we found are, in the order of the cases in the theorem, 2, 3, 6, 4 and 3.
In most cases, $Q$ depends on coefficients $a$ which are subject to
(solvable!) systems of polynomial equations. We list those matrices on the
web-page \cite{RWwebpage}, in the form of a Singular-executable text file.
This page also provides a few small Singular routines to perform the necessary
checks: extraction of the matrix factorisation conditions (bilinear equations
on the $a$), computation of the quantum dimensions, computation of
the Gr\"obner basis for the ideal in (\ref{idealmemb}).
For the sake of completeness, and in order to give an impression of
the complexity, the matrices and the polynomial equations are also
reproduced in the appendix of the present paper.

In all of the five cases, the orbifold equivalence satisfies
$q_L(Q)q_R(Q) \neq \pm 1$, hence $A = Q^\dagger\star Q$ is not
similar to the identity defect: these are ``true orbifold
equivalences'', not ``mere equivalences in the bicategory
${\cal LG}$''. $\quad\square$
\vspace{.2cm}

The web-page mentioned above also presents direct orbifold equivalences
between $D_7$ and $E_6$, between $D_{10}$ and $E_7$,
and between $D_{16}$ and $E_8$.
That these simple singularities are orbifold equivalent follows already from
the $A$-$D$ and $A$-$E$ results in \cite{CR3,CRCR}, what makes the direct
$D$-$E$ defects noteworthy is that they have at most rank 3. (The smallest
orbifold equivalence between $E_8$ and $A_{29}$ is of rank 4.)
\vspace{.3cm}

Together with the straightforward $E_{14}$-$Q_{10}$ orbifold equivalence
mentioned in section 2, Theorem 4.1 exhausts all orbifold equivalences
among the (quasi-homogeneous) exceptional unimodal Arnold singularities:
no other pairs with equal central charge exist among those fourteen
potentials. The orbifold equivalent pairs are precisely the pairs
that display ``strange duality'' (Dolgachev and Gabrielov numbers
are interchanged), see e.g.\ \cite{Saito}. 

Among the 14 exceptional bimodal Arnold singularities, only $Q_{17}$
and $W_{17}$ have the same central charge (namely $\hat c = {6\over5}$.);
we have not yet found an orbifold equivalence between them (nor between
$Q_{17}$ or $W_{17}$ and the pair in item (5) above). 

\vspace{.1cm}
It might be worth mentioning that the arguments one can use to treat 
the $E_{14}$-$Q_{10}$ case -- i.e.\ Theorem 2.2 (c) -- also show that
orbifold equivalence does not respect the modality of a singularity:

\noindent
The exceptional unimodal Arnold singularity $Q_{12}$ with 
$V_{Q_{12}}(x) = x_1^5+x_1\,x_2^2 +x_3^3$ is orbifold equivalent to the
exceptional bimodal Arnold singularity $E_{18}$ given by 
$V_{E_{18}}(y) = y_1^{10}+ y_2^3  +y_3^2$: the former is $D_6 \times A_2$, the
latter $A_9 \times A_2$, and $D_6 \oeq\ A_9$ due to the results of \cite{CR3}.
By the same method, one can relate other exceptional Arnold
singularities to sums of simple singularities; among the examples
involving bimodal singularities are $Q_{16} \oeq A_{13} \times A_{2}$ 
and $U_{16} \oeq E_8 \times A_2 \oeq A_5 \times A_4$.

\vspace{.3cm}
A number of more or less expected orbifold equivalences, including infinite
series, can be established via transformations of variables:

\vspace{.2cm} 
\noindent
\textbf{Lemma 4.2: } Assume $Q(x,y)$ is an orbifold equivalence between
$V_1(x)$ and $V_2(y)$, and assume that $y\mapsto y'$ is an invertible,
weight-preserving transformation of variables. Then $Q(x,y')$ is an
orbifold equivalence between $V_1(x)$ and $V_2(y')$ if the weights
$|y_i|$ are pairwise different, or if $V_2(y) \in \maxm^3$. 

\vspace{.2cm} 
\noindent
\textsl{Proof:} First, focus on the variable transformation itself: We can assume
wlog. that the $y_1,\ldots,y_m$ are labeled by increasing weight, $y_1$ having the
lowest weight. Then the transformation can be written as
$y_j \mapsto y_j' = f_j(y) + \sum_{k\in I_j} A_{jk} y_k$ where $A_{jk}\in \C$, where
$I_j = \{ k\, :\, |y_k| = |y_j|\}$ and where $f_j$ depends only on those $y_l$ with
$|y_l| < |y_j|$. As $y\mapsto y'$ preserves weights, $f_j$ has no linear terms.
The Jacobian ${\cal J}$ of the transformation is lower block-diagonal and
$\det({\cal J}) = \det(A)$, a non-zero constant. 

Since $Q' := Q(x,y')$ is obviously a matrix factorisation of
$V_1(x)-V_2(y')$, we only need to study the quantum dimensions of $Q'$.
The relation $q_R(Q') = \det(A)\,q_R(Q)$ results immediately from making a
substitution of integration variables in the formula for the right quantum
dimension.

The left quantum dimension of $Q'$ can be expressed as a Kapustin-Li correlator
(in the $(V_1,Q_1)$ model) of the fermions $F'_j = \partial_{y'_j}Q'|_{y'=0} 
= \sum_{l=1}^m \frac{\partial y_l}{\partial y'j}|_{y=0}\,F_l$.
Here, we have already exploited
$y'=0 \Leftrightarrow y=0$ to simplify, but the summation over $l$ might
still lead to linear combinations which are difficult to control. The extra
assumptions on $V_2(y)$ avoid this: If all $|y_j|$ are pairwise different,
then $\frac{\partial y_l}{\partial y'j}|_{y=0} = b_j\,\delta_{j,l}$ for some
non-zero constants $b_j$. If $V_2$ starts at order 3 or higher, the
$F_j$ anti-commute with each other inside the correlator: adapting the
proof of Prop.\ 3.3, one finds 
$$ 0 = -  \partial_{y_{j_1}} \partial_{y_{j_2}} V_2|_{y=0} = 
\{F_{j_1},F_{j_2}\} + \{Q_1(x), \partial_{y_{j_1}} \partial_{y_{j_2}} Q(x,y)|_{y=0}\}\ ,
$$
and the last term vanishes in the $Q_1$-cohomology, therefore does not contribute
to the Kapustin-Li correlator. Hence, the correlator is totally anti-symmetric
in the $F_j$, and the linear combination of correlators making up the left
quantum dimension is simply $q_L(Q') = \det(A)^{-1}\,q_L(Q)$. 
$\quad \square$ 

\vspace{.3cm}
Applying this lemma to the identity defect of $V_1(x)-V_1(y)$, one can
establish orbifold equivalences e.g.\ in the following cases: 

\begin{itemize}

\item[(1)]
So-called ``auto-equivalences'' of unimodal Arnold singularities: different
descriptions of the same singularity exist for
$U_{12},\;Q_{12},\;W_{12},\;W_{13},\;Z_{13}$ and
$E_{14}$. The assumptions on the variable weights resp.\ structure of $V_2$
made in Lemma 4.2 hold for all these cases. These orbifold equivalences were
already
discussed in \cite{NRC2}, and although the concrete formulas given there contain
errors, the general structure ($Q$ being a nested tensor product of rank 4)
coincides with what one obtains from the identity defect upon a weight-preserving
transformation of variables.

The auto-equivalence between $V_{Q_{17}^T}(x) = x_1^3x_2+x_2^5 x_3 +x_3^2$
and $V(y) = y_1^3 y_2 + y_2^{10} +y_3^2$ is another such example, involving
a bimodal Arnold singularity.

\item[(2)] Equivalences between quasi-homogeneous polynomials of Fermat, chain
  and loop (or cycle) type at $\hat c <1$: 

  $V_{A_{2n-1}}(x) = x_1^{2n} + x_2^2\ $ and $\ V_{D^T_{n+1}}(y) = y_1^n y_2 + y_2^2$

  $V_{L_n}(x) = x_1^nx_2 + x_1x_2^2\ $ and $\ V_{D_{2n}}(y) = y_1^{2n-1} + y_1 y_2^2$  

  $V_{C_n}(x) = x_1^2x_2 + x_2^nx_3+x_3^2\ $ and $\ V_{D_{2n+1}}(y) = y_1^{2n} + y_1 y_2^2 + y_3^2$

  with $n\geq 2$ in all three pairs. Explicit orbifold equivalences for
  $A$-$D^T$ were already given in \cite{RCThesis}. 

\item[(3)] Cases involving non-trivial marginal bulk deformations, e.g.\ 

  at central charge $\hat c = \frac{10}{9}$, one finds an orbifold
  equivalence between the product $A_8 \times A_2$ of 
  simple singularities, $V_{(A_8\times A_2)} (x) =  x_1^9 +x_2^3$,
  and special deformations of $Z_{13}^T$, given by
  $V_{Z_{13}^T}(y) = y_1^6 y_2 +y_2^3 + \mu_2\, y_1^3 y_2^2$,
  if $\mu_2 = \pm\sqrt{3}$;

at central charge $\hat c = \frac{8}{7}$,  the two deformed singularities
$V_{E_{19}^T}(x) = x_1^3x_2+ x_2^7 + \mu_1\,x_1 x_2^5$ and
$V_2(y) = y_1 y_2^5 +y_1^3 y_2 + \mu_2\, y_1^2 y_2^3$ are
orbifold equivalent as long as the two deformation parameters are
related by $\mu_1 = \mu\,({1\over3}\mu_2^2 - 1)$ with
$3\mu^3 = -\mu_2({2\over9}\mu_2^2 - 1)$. 

\end{itemize}

\vspace{.2cm} 
Lemma 4.2 can also be used to prove an orbifold equivalence one would
expect on geometric grounds: The elliptic curve is  described by 
$V_{\lambda}(x) = - x_2^2 x_3 + x_1 (x_1-x_3) (x_1-\lambda x_3)$,
where $\lambda$ is a complex parameter with $\lambda \neq 0,1$,  
and two such curves $V_\lambda$ and $V_{\lambda'}$ describe birationally 
equivalent tori if and only if
\be\la{latildevalues}
\lambda' \in \{\, \lambda,\; 1-\lambda,\;1/\lambda, \;1/(1-\lambda),
\;(\lambda-1)/\lambda, \;\lambda/(\lambda-1)\,\}\ . 
\ee
One can apply a weight-preserving variable transformation to bring
$V_\lambda$ into an alternative form
$V_e(y) = - y_2^2 y_3 + (y_1-e_1\,y_3)(y_1-e_2\,y_3)(y_1-e_3\,y_3)$. 
The parameters of the two forms are related by
$\lambda = (e_3-e_1)/(e_2-e_1)$, and the
six different $\lambda'$-values in (\ref{latildevalues}) arise from 
permuting the $e_i$, which of course leaves $V_e$ unchanged; 
thus we find $V_{\lambda'} \oeq V_\lambda$.

\vspace{.2cm}
Since, in all the examples listed after Lemma 4.2, we start from the
identity defect,
the orbifold equivalence resulting from the transformation of variables
automatically satisfies $q_L(Q')\,q_R(Q') = 1$, so it is likely that
they are ``mere equivalences'' in the bicategory ${\cal LG}$.
(One way to verify this would be to compute and analyse the fusion product
$(Q')^\dagger \star Q'$.) But Lemma 4.2 could also be applied to the
orbifold equivalence between $D_{n+1}$ and $A_{2n-1}$, say, to produce 
a defect with $q_L(Q')\,q_R(Q') = 2$ between $D_{n+1}$ and $D^T_{n+1}$. 

Furthermore, the potentials of type $D^T_n$, $C_n$ and $L_n$ listed in
item (2) appear as separate entries in lists of quasi-homogeneous polynomials
\cite{KrSk,HertlingKurbel}, but not in lists of singularities
(where more general types of transformations of variables are allowed to
identify two singularities).  The orbifold equivalences given in item (2)
of Lemma 4.2 may not be surprising, but it is not clear to us whether
there are abstract theorems guaranteeing that polynomials which are
equivalent as singularities are (orbifold) equivalent in  ${\cal LG}$. 

\vspace{.2cm}
A first edition of an ``oeq catalogue'', i.e.\ a list of polynomials
sorted into orbifold equivalence classes based on the results of
\cite{CR3,CRCR} and our new findings, is available at the web-page
\cite{RWwebpage}.

\section{Open problems and conjectures}\label{sec:examples}

Ultimately, one would like to find a simple (combinatorial or
number-theoretic) criterion that
allows to read off directly from the potentials $V_1,\,V_2$
whether they are orbifold equivalent or not -- instead of taking
a detour via constructing an explicit orbifold equivalence $Q$. 

Having invested quite a lot of effort into finding such matrices
$Q$, the
authors sincerely hope that such a criterion involves conditions
beyond the ones listed in Prop.\ 2.3. 

And there are indeed reasons to believe that $\hat c(V_1) = \hat c(V_2)$
alone is insufficient for $V_1 \oeq V_2$: 

One line of arguments concerns marginal deformations:
Le $V$ be a potential which admits a marginal deformation, i.e.\ there
is a quasi-homogeneous element $\phi \in \mathrm{Jac}(V)$ of weight 2.
(The Fermat elliptic curve $V(x) = \sum_{i=1}^3 x_i^3$ with
$\phi = x_1x_2x_3$ is an example.) Set $V_1(x) = V(x) + \mu\,\phi(x)$,
where $\mu\in \C$ is a deformation parameter, and $V_2(y) = V(y)$; we
have $\hat c(V_1) = \hat c(V_2)$.

The examples of orbifold equivalences listed at the end of section 4,
involving $Z_{13}^T$, $E_{19}^T$ and the geometrically equivalent tori
(\ref{latildevalues}), already suggest that a given method of constructing
a defect $Q$ between $V+\mu \phi$ and $V$ might lead to an orbifold
equivalence for a discrete set of $\mu$-values only. 

In general, let $Q(x,y;\mu)$ be rank $N$ matrix factorisation of $V_1(x)-V_2(y)$
and assume its $\mu$-derivative exists in a neighbourhood of $\mu=0$. Then
the bosonic morphism $\Phi := \phi(x)\,\mathbf{1}_{2N} = \{Q,\partial_\mu Q\}$
is zero in the cohomology of $Q$, and $\mathbf{1}\otimes\Phi$ is zero in the
cohomology of $Q^\dagger \hat\otimes Q$. This should imply that
$\Phi$ is absent from $\mathrm{End}(A)$ for $A = Q^\dagger \star Q$,
which in turn makes it unlikely that there is a projection 
from $\mathrm{End}(A)$ to $\mathrm{Jac}(V)$ -- but the latter has to be
the case \cite{CR3} if $Q$ is an orbifold equivalence.

A more direct proof that $c(V_1)=\hat c(V_2)$ does not guarantee
orbifold equivalence might result from incompatibility of the ``weight split
lists'' $S_\tau$ occurring in the weight split criterion from subsection 3.2.

\vspace{.2cm} 
We conjecture that orbifold equivalences $Q$ have trivial fermionic cohomology.

This is true in every concrete case for which we have computed $H^1_{Q}$,
and the conjecture is backed up by the following observation: If a matrix
factorisation $Q$ of $W$ has a non-trivial fermion $\psi \in H^1_{Q}$, one
can form the cone
$$ 
C_\psi(\lambda) =  \left(
  \begin{array}{cc}
            Q  & \lambda\psi \\
            0 &  Q   
\end{array}
 \right)
$$
which is again a matrix factorisation of $W$ for any $\lambda \in \C$.
The upper triangular form implies that
$q_L(C_\psi(\lambda)) = q_L(Q\,\oplus\,Q) =2 q_L(Q)$ 
for any value of $\lambda$ (likewise for the right quantum dimension).
In general, however, cones $C_\psi(\lambda)$ with $\lambda \neq 0$ are
not equivalent (related by similarity transformations) to the direct
sum $Q\,\oplus\,Q$, so one would not expect the quantum dimensions to
always coincide. 

A related question (related due to the role of fermions in deformations
of matrix factorisations) is whether there can be moduli spaces of
orbifold equivalences between two fixed potentials, or whether the
equations only ever admit a discrete set of solutions. Our computations
point towards the latter, but we have no proof. 

The bicategory setting might provide a better language in which
to tackle these general questions.

\vspace{.2cm}
We hope that the orbifold equivalences presented here prove fruitful in
singularity theory, and in other areas related to matrix factorisations
by well-established equivalences of categories, but one should also
explore applications of orbifold equivalence in string theory, or
in the context of mirror symmetry.

E.g., orbifold equivalences between Arnold singularities may also imply relations
between $N=2$ supersymmetric gauge theories in 4 dimensions ``engineered''
from these singularities, see \cite{CecDelZ,DelZ}. 

It is reasonable to expect that potentials related by
``Berglund-H\"ubsch-Krawitz'' duality \cite{BergHu,Kraw}
are orbifold equivalent ($E_{13} \oeq Z_{11}$ is one example),
and this duality is one approach to constructing mirror manifolds.

\vspace{.2cm}
One might also explore whether some of the orbifold equivalences of
Landau-Ginzburg potentials can be ``lifted'' to relations of the
conformal field theories associated with them. In particular, 
the question whether there is a CFT analogue to $A_5 \oeq A_2 \times A_2$,
perhaps in terms of an orbifold construction, should be accessible because
the central charge is that of a theory of
two free bosons. The $A_2$-model is assocatiated with a free boson
compactified on a circle, see \cite{LVW} and references therein; we
do not know whether a similar statement can be made for $A_5$.

\vspace{.3cm}
Let us add some speculative comments on orbifold equivalence and
entanglement. In quantum physics, entanglement refers to the phenomenon
that a physical system comprised from two subsystems (like an
electron-positron pair) can be in a state such that observations made
on one subsystem immediately determine properties of the second subsystem
no matter how great the separation between the two. This behaviour has
no analogue in classical physics.

Already the general consequences implied by an
orbifold equivalence $V_1 \oeq V_2$ -- e.g.\ the relation
(\ref{hmfmodequiv}) between categories, or more directly the one between
correlators in the Landau-Ginzburg models associated with $V_1(x)$ and $V_2(y)$ --
are strongly reminiscent of entanglement.

Closer to the level of concrete formulas, one notices that quantum states 
displaying entanglement are formed from states describing the subsystems in a
manner that resembles the mixing of $x$- and $y$-variables implied by Prop.\ 3.3.

Indeed, we expect that the quantum dimensions of a defect can be
related to a suitably defined entanglement entropy in Landau-Ginzburg models.

If this can be made manifest and the ``symmetries'' discussed here can
ultimately be traced back to quantum entanglement, perhaps
``entanglement equivalence'' might be a more appropriate term than
``orbifold equivalence''.

\subsection*{Acknowledgements}

We are indebted to N.~Carqueville for introducing us to the problem,
for many valuable discussions, and for an early version of the Singular
code to compute quantum dimensions. We also thank I.~Brunner, D.~Murfet and
I.~Runkel for useful conversations and comments.

The work of P.W.\ was partially supported by an STFC studentship.

\begin{appendix}
\section{Appendix:\ explicit defects}

For the sake of completeness, we collect the orbifold equivalences that
can serve to prove Theorem 4.1. The Singular-executable formats given on
the web-page \cite{RWwebpage} should be of more practical use. 

To save writing zeroes, we list matrices $E$ and $J$ only.
$Q$ is constructed from them as in (\ref{QEJ}). For fear of producing
typos, we have largely refrained from attempts at simplifying the Singular output
(except for the very easy case $A_5 \oeq A_2 \times A_2$).
The matrices spelled out in the
following are the simplest ones we could find: what results from our
Singular algorithm typically contains many more coefficients $a_{rs,\vec p}$,
and we have chosen explicit values for some of them. 

The orbifold equivalences are listed in the order they appear in Theorem 4.1.

\vspace{.5cm}
\noindent
\textbf{(1)} A rank 2 orbifold equivalence between $A_5$ and $A_2 \times A_2$:

\small{
$$
  E= \left(
  \begin{array}{cc}
    x_1^2 -a_1( y_1 + y_2)
    & \ x_2 + a_2 x_1 (y_1- y_2)  \\
    x_2 - a_2 x_1 (y_1- y_2) 
    & \ - x_1^4 -64 a_1^8 y_2^2 +16 a_1^5 y_1 y_2 
    -a_1 x_1^2 (y_1 + y_2) -4 a_1^2 y_1^2 
   \end{array}
 \right)
$$

$$
  J= \left(
  \begin{array}{cc}
    x_1^4 + 64 a_1^8 y_2^2 -16 a_1^5 y_1 y_2 
    + a_1 x_1^2  (y_1 + y_2)   +4 a_1^2 y_1^2 
    & \ x_2 + a_2 x_1 (y_1  - y_2) 
    \\   
    x_2 - a_2 x_1 (y_1 - y_2) 
    & \  - x_1^2+ a_1 (y_1 + y_2)        
   \end{array}
 \right)
$$
}
\normalsize{
where the coefficients have to satisfy 
$$ a_2^2 =  3 a_1^2 \quad \mathrm{and}\quad  a_1^3 = {\textstyle {1\over 4}}\ .
$$
The quantum dimensions of $Q$ are $q_L(Q) = -2a_1a_2$
and $q_R(Q) = -{4\over3} a_2$. Since their product is 2, this is a ``true
orbifold equivalence'', not an ordinary equivalence in the bi-category of
Landau-Ginzburg potentials. On the other hand, since 2 is contained in
any cyclotomic field, a group action might be the source of this orbifold
equivalence.
}

\vspace{.5cm}

\normalsize{
\noindent  
\textbf{(2)} An ugly rank 3 orbifold equivalence between $E_{13}$ and $Z_{11}$:   }

\noindent The matrix elements $E_{rs}$ and $J_{rs}$ are given by 

\vspace{.2cm}   
\small{
\noindent $
E_{11}=-x_1^2-y_1 a_2
$ \newline\noindent $
E_{12}=-x_1 y_2 a_3-x_1 y_2 a_4+x_2
$ \newline\noindent $
E_{13}=-y_2 a_4
$ \newline\noindent $
E_{21}=x_1 y_2 a_3-x_1 y_2 a_5+x_2
$ \newline\noindent $
E_{22}=y_2^2 a_1^2 a_4^2+y_2^2 a_1 a_3 a_4+y_2^2 a_1 a_4^2+y_2^2 a_1 a_4 a_5-x_1^3 a_1+y_2^2 a_3^2+y_2^2 a_3 a_4-y_2^2 a_3 a_5+y_2^2 a_5^2+x_1 y_1 a_1 a_2-x_1 y_1 a_1 a_6+x_1^3-y_2^2 a_7
$ \newline\noindent $
E_{23}=x_1^2+y_1 a_6
$ \newline\noindent $
E_{31}=-x_1^3 a_1-x_1^3-y_2^2 a_7
$ \newline\noindent $
E_{32}=-x_1^2 y_2 a_1^2 a_4+2 y_1 y_2 a_1^2 a_4 a_2-x_1^2 y_2 a_1 a_3-x_1^2 y_2 a_1 a_4-x_1^2 y_2 a_1 a_5+y_1 y_2 a_1 a_3 a_2+2 y_1 y_2 a_1 a_4 a_2+y_1 y_2 a_1 a_5 a_2-y_1 y_2 a_1 a_4 a_6-x_1^2 y_2 a_3+x_1^2 y_2 a_5+y_1 y_2 a_3 a_2-y_1 y_2 a_4 a_6-y_1 y_2 a_5 a_6+x_1 x_2
$ \newline\noindent $
E_{33}=x_1 y_2 a_5+x_2
$

\vspace{.2cm}
\noindent $
J_{11}=-x_1 y_2^3 a_1^2 a_4^2 a_5-x_1^4 y_2 a_1^2 a_4-x_1 y_2^3 a_1 a_3 a_4 a_5-x_1 y_2^3 a_1 a_4^2 a_5-x_1 y_2^3 a_1 a_4 a_5^2+2 x_1^2 y_1 y_2 a_1^2 a_4 a_2-x_1^2 y_1 y_2 a_1^2 a_4 a_6+2 y_1^2 y_2 a_1^2 a_4 a_2 a_6-x_1^4 y_2 a_1 a_3-x_1^4 y_2 a_1 a_4-x_2 y_2^2 a_1^2 a_4^2-x_1 y_2^3 a_3^2 a_5-x_1 y_2^3 a_3 a_4 a_5+x_1 y_2^3 a_3 a_5^2-x_1 y_2^3 a_5^3+x_1^2 y_1 y_2 a_1 a_3 a_2+2 x_1^2 y_1 y_2 a_1 a_4 a_2-x_1^2 y_1 y_2 a_1 a_3 a_6-2 x_1^2 y_1 y_2 a_1 a_4 a_6+y_1^2 y_2 a_1 a_3 a_2 a_6+2 y_1^2 y_2 a_1 a_4 a_2 a_6+y_1^2 y_2 a_1 a_5 a_2 a_6-y_1^2 y_2 a_1 a_4 a_6^2-x_1^4 y_2 a_3-x_2 y_2^2 a_1 a_3 a_4-x_2 y_2^2 a_1 a_4^2-x_2 y_2^2 a_1 a_4 a_5+x_1^2 y_1 y_2 a_3 a_2-x_1^2 y_1 y_2 a_3 a_6-x_1^2 y_1 y_2 a_4 a_6+y_1^2 y_2 a_3 a_2 a_6-y_1^2 y_2 a_4 a_6^2-y_1^2 y_2 a_5 a_6^2+x_1 y_2^3 a_5 a_7+x_1^3 x_2 a_1-x_2 y_2^2 a_3^2-x_2 y_2^2 a_3 a_4+x_2 y_2^2 a_3 a_5-x_2 y_2^2 a_5^2-x_1 x_2 y_1 a_1 a_2+x_1 x_2 y_1 a_1 a_6+x_1 x_2 y_1 a_6+x_2 y_2^2 a_7 
$ \newline\noindent $
J_{12}=-x_1^2 y_2^2 a_1^2 a_4^2-y_1 y_2^2 a_1^2 a_4^2 a_2-x_1^2 y_2^2 a_1 a_3 a_4-x_1^2 y_2^2 a_1 a_4^2-x_1^2 y_2^2 a_1 a_4 a_5-y_1 y_2^2 a_1 a_3 a_4 a_2-y_1 y_2^2 a_1 a_4^2 a_2-y_1 y_2^2 a_1 a_4 a_5 a_2-x_1^2 y_2^2 a_3 a_4-x_1^2 y_2^2 a_3 a_5-y_1 y_2^2 a_3^2 a_2-y_1 y_2^2 a_3 a_4 a_2+y_1 y_2^2 a_3 a_5 a_2-y_1 y_2^2 a_5^2 a_2-x_1 y_1^2 a_1 a_2^2+x_1 y_1^2 a_1 a_2 a_6+x_1 y_1^2 a_2 a_6+y_1 y_2^2 a_2 a_7-y_1 y_2^2 a_6 a_7-x_1 x_2 y_2 a_3+x_1 x_2 y_2 a_5+x_2^2
$ \newline\noindent $
J_{13}=-y_2^3 a_1^2 a_4^3-y_2^3 a_1 a_3 a_4^2-y_2^3 a_1 a_4^3-y_2^3 a_1 a_4^2 a_5+x_1^3 y_2 a_1 a_4-y_2^3 a_3^2 a_4-y_2^3 a_3 a_4^2+y_2^3 a_3 a_4 a_5-y_2^3 a_4 a_5^2-x_1 y_1 y_2 a_1 a_4 a_2+x_1 y_1 y_2 a_1 a_4 a_6+x_1^3 y_2 a_3+x_1 y_1 y_2 a_3 a_6+x_1 y_1 y_2 a_4 a_6+y_2^3 a_4 a_7-x_1^2 x_2-x_2 y_1 a_6
$ \newline\noindent $
J_{21}=-3 y_1 y_2^2 a_1^2 a_4^2 a_2-2 y_1 y_2^2 a_1 a_3 a_4 a_2-3 y_1 y_2^2 a_1 a_4^2 a_2-2 y_1 y_2^2 a_1 a_4 a_5 a_2+y_1 y_2^2 a_1 a_4^2 a_6+x_1^5 a_1+x_1^2 y_2^2 a_3 a_5-x_1^2 y_2^2 a_5^2-y_1 y_2^2 a_3^2 a_2-2 y_1 y_2^2 a_3 a_4 a_2+y_1 y_2^2 a_3 a_5 a_2-y_1 y_2^2 a_5^2 a_2-x_1 y_1^2 a_1 a_2^2+x_1^3 y_1 a_1 a_6+y_1 y_2^2 a_4^2 a_6+y_1 y_2^2 a_4 a_5 a_6+x_1 y_1^2 a_1 a_2 a_6+x_1^5+x_1^3 y_1 a_6+x_1 y_1^2 a_2 a_6+x_1^2 y_2^2 a_7+y_1 y_2^2 a_2 a_7+x_1 x_2 y_2 a_3+x_2^2
$ \newline\noindent $
J_{22}=x_1^3 y_2 a_1 a_4+x_1^3 y_2 a_4+x_1^3 y_2 a_5+x_1 y_1 y_2 a_5 a_2+y_2^3 a_4 a_7+x_1^2 x_2+x_2 y_1 a_2
$ \newline\noindent $
J_{23}=x_1 y_2^2 a_3 a_4-x_1 y_2^2 a_4 a_5-x_1^4-x_1^2 y_1 a_2-x_1^2 y_1 a_6-y_1^2 a_2 a_6+x_2 y_2 a_4
$ \newline\noindent $
J_{31}=-x_1^3 y_2^2 a_1^3 a_4^2-2 x_1^3 y_2^2 a_1^2 a_4^2-2 x_1^3 y_2^2 a_1^2 a_4 a_5-2 x_1 y_1 y_2^2 a_1^2 a_3 a_4 a_2+3 x_1 y_1 y_2^2 a_1^2 a_4^2 a_2+2 x_1 y_1 y_2^2 a_1^2 a_4 a_5 a_2-y_2^4 a_1^2 a_4^2 a_7+x_1^6 a_1^2-x_1^3 y_2^2 a_1 a_3 a_4-x_1^3 y_2^2 a_1 a_4^2+x_1^3 y_2^2 a_1 a_3 a_5-2 x_1^3 y_2^2 a_1 a_4 a_5-2 x_1^3 y_2^2 a_1 a_5^2-x_1^4 y_1 a_1^2 a_2-x_1 y_1 y_2^2 a_1 a_3^2 a_2+3 x_1 y_1 y_2^2 a_1 a_4^2 a_2+4 x_1 y_1 y_2^2 a_1 a_4 a_5 a_2+x_1 y_1 y_2^2 a_1 a_5^2 a_2+x_1^4 y_1 a_1^2 a_6+x_1 y_1 y_2^2 a_1 a_3 a_4 a_6-x_1 y_1 y_2^2 a_1 a_4^2 a_6-x_1 y_1 y_2^2 a_1 a_4 a_5 a_6-y_2^4 a_1 a_3 a_4 a_7-y_2^4 a_1 a_4^2 a_7-y_2^4 a_1 a_4 a_5 a_7+x_1^2 x_2 y_2 a_1^2 a_4-x_1^3 y_2^2 a_3 a_4-x_1^3 y_2^2 a_3 a_5-x_1^4 y_1 a_1 a_2-2 x_2 y_1 y_2 a_1^2 a_4 a_2+2 x_1 y_1 y_2^2 a_3 a_4 a_2+x_1 y_1 y_2^2 a_5^2 a_2+x_1^2 y_1^2 a_1 a_2^2+x_1^4 y_1 a_1 a_6+x_1 y_1 y_2^2 a_3 a_4 a_6-x_1 y_1 y_2^2 a_4^2 a_6+x_1 y_1 y_2^2 a_3 a_5 a_6-2 x_1 y_1 y_2^2 a_4 a_5 a_6-x_1 y_1 y_2^2 a_5^2 a_6-x_1^2 y_1^2 a_1 a_2 a_6+2 x_1^3 y_2^2 a_1 a_7-y_2^4 a_3^2 a_7-y_2^4 a_3 a_4 a_7+y_2^4 a_3 a_5 a_7-y_2^4 a_5^2 a_7-x_1 y_1 y_2^2 a_1 a_2 a_7+x_1 y_1 y_2^2 a_1 a_6 a_7-x_1^6+x_1^2 x_2 y_2 a_1 a_3+x_1^2 x_2 y_2 a_1 a_4+x_1^2 x_2 y_2 a_1 a_5-x_2 y_1 y_2 a_1 a_3 a_2-2 x_2 y_1 y_2 a_1 a_4 a_2-x_2 y_1 y_2 a_1 a_5 a_2+x_2 y_1 y_2 a_1 a_4 a_6-x_1^2 y_1^2 a_2 a_6-x_1 y_1 y_2^2 a_2 a_7+x_1 y_1 y_2^2 a_6 a_7+y_2^4 a_7^2-x_2 y_1 y_2 a_3 a_2+x_2 y_1 y_2 a_4 a_6+x_2 y_1 y_2 a_5 a_6-x_1 x_2^2
$ \newline\noindent $
J_{32}=x_1^4 y_2 a_1^2 a_4-x_1^2 y_1 y_2 a_1^2 a_4 a_2-2 y_1^2 y_2 a_1^2 a_4 a_2^2+x_1^4 y_2 a_1 a_5-x_1^2 y_1 y_2 a_1 a_4 a_2-y_1^2 y_2 a_1 a_3 a_2^2-2 y_1^2 y_2 a_1 a_4 a_2^2-y_1^2 y_2 a_1 a_5 a_2^2+x_1^2 y_1 y_2 a_1 a_4 a_6+y_1^2 y_2 a_1 a_4 a_2 a_6-x_1^4 y_2 a_4-x_1^4 y_2 a_5-x_1^2 y_1 y_2 a_5 a_2-y_1^2 y_2 a_3 a_2^2+x_1^2 y_1 y_2 a_4 a_6+x_1^2 y_1 y_2 a_5 a_6+y_1^2 y_2 a_4 a_2 a_6+y_1^2 y_2 a_5 a_2 a_6-x_1 y_2^3 a_3 a_7-x_1 y_2^3 a_4 a_7+x_1^3 x_2 a_1-x_1 x_2 y_1 a_2+x_2 y_2^2 a_7
$ \newline\noindent $
J_{33}=x_1^2 y_2^2 a_1^2 a_4^2-2 y_1 y_2^2 a_1^2 a_4^2 a_2+x_1^2 y_2^2 a_1 a_3 a_4+x_1^2 y_2^2 a_1 a_4^2+x_1^2 y_2^2 a_1 a_4 a_5-y_1 y_2^2 a_1 a_3 a_4 a_2-2 y_1 y_2^2 a_1 a_4^2 a_2-y_1 y_2^2 a_1 a_4 a_5 a_2+y_1 y_2^2 a_1 a_4^2 a_6-x_1^5 a_1+x_1^2 y_2^2 a_4 a_5+x_1^2 y_2^2 a_5^2-y_1 y_2^2 a_3 a_4 a_2-x_1^3 y_1 a_1 a_6+y_1 y_2^2 a_4^2 a_6+y_1 y_2^2 a_4 a_5 a_6+x_1^5+x_1^3 y_1 a_2+x_1 y_1^2 a_2 a_6-x_1^2 y_2^2 a_7-y_1 y_2^2 a_6 a_7-x_1 x_2 y_2 a_4-x_1 x_2 y_2 a_5+x_2^2$
}

\vspace{.3cm}
\normalsize{
  The seven coefficients $a_1,\ldots,a_7$ are subject to matrix factorisation
  conditions which take the form of twelve algebraic equations $f_{\alpha}(a)=0$ with
}
\vspace{.1cm}

\small{
\noindent 
$f_1 =
     -(1/3) a_1 a_3^2 a_4 a_6-(1/3) a_1 a_3 a_4^2 a_6+(1/3) a_1 a_3 a_4 a_5 a_6+(2/3) a_1 a_4^2 a_5 a_6+(2/3) a_1 a_4 a_5^2 a_6-(2/3) a_2 a_3^3-(1/3) a_2 a_3^2 a_4+2 a_2 a_3^2 a_5+(2/3) a_2 a_3 a_4 a_5-2 a_2 a_3 a_5^2+(4/3) a_2 a_5^3-(1/3) a_3^2 a_4 a_6-(1/3) a_3 a_4^2 a_6+(2/3) a_4^2 a_5 a_6+(4/3) a_4 a_5^2 a_6+(2/3) a_2 a_3 a_7-(4/3) a_2 a_5 a_7-(5/3) a_3 a_6 a_7+(1/3) a_5 a_6 a_7
$
\newline\noindent 
$f_2 =
2 a_1^2 a_3 a_4^2 a_6-4 a_1^2 a_4^2 a_5 a_6+2 a_1 a_3 a_4^2 a_6-4 a_1 a_4^2 a_5 a_6-2 a_2 a_3^3+6 a_2 a_3^2 a_5-6 a_2 a_3 a_5^2+4 a_2 a_5^3-2 a_3 a_4 a_5 a_6+4 a_4 a_5^2 a_6+2 a_2 a_3 a_7-4 a_2 a_5 a_7-6 a_3 a_6 a_7
$
\newline\noindent 
$f_3 =
-7 a_1^3 a_2 a_3 a_4^2 a_6+11 a_1^3 a_2 a_4^2 a_5 a_6-4 a_1^2 a_2 a_3^2 a_4 a_6-12 a_1^2 a_2 a_3 a_4^2 a_6-2 a_1^2 a_2 a_3 a_4 a_5 a_6+15 a_1^2 a_2 a_4^2 a_5 a_6+14 a_1^2 a_2 a_4 a_5^2 a_6-a_1^2 a_4^3 a_6^2-a_1 a_2 a_3^3 a_6-9 a_1 a_2 a_3^2 a_4 a_6-8 a_1 a_2 a_3 a_4^2 a_6+4 a_1 a_2 a_3 a_4 a_5 a_6+4 a_1 a_2 a_4^2 a_5 a_6+7 a_1 a_2 a_4 a_5^2 a_6+5 a_1 a_2 a_5^3 a_6-2 a_1 a_4^3 a_6^2-2 a_1 a_4^2 a_5 a_6^2-12 a_1^2 a_2 a_4 a_6 a_7-a_2^2 a_3^3-2 a_2^2 a_3^2 a_4+a_2^2 a_3^2 a_5-a_2^2 a_3 a_5^2-a_2 a_3^3 a_6-5 a_2 a_3^2 a_4 a_6-3 a_2 a_3 a_4^2 a_6+a_2 a_3^2 a_5 a_6+a_2 a_3 a_4 a_5 a_6-a_2 a_3 a_5^2 a_6+a_2 a_5^3 a_6-a_4^3 a_6^2-2 a_4^2 a_5 a_6^2-a_4 a_5^2 a_6^2-6 a_1 a_2 a_3 a_6 a_7-7 a_1 a_2 a_4 a_6 a_7-6 a_1 a_2 a_5 a_6 a_7+a_1 a_4 a_6^2 a_7+a_2^2 a_3 a_7-a_2 a_3 a_6 a_7-a_2 a_5 a_6 a_7+a_4 a_6^2 a_7+a_5 a_6^2 a_7
$
\newline\noindent 
$f_4 =
      -(5/2) a_1^3 a_3 a_4^2+2 a_1^3 a_4^2 a_5-2 a_1^2 a_3^2 a_4-(9/2) a_1^2 a_3 a_4^2+2 a_1^2 a_3 a_4 a_5+3 a_1^2 a_4^2 a_5-2 a_1^2 a_4 a_5^2-(3/2) a_1 a_3^3-4 a_1 a_3^2 a_4-2 a_1 a_3 a_4^2+(9/2) a_1 a_3^2 a_5+(11/2) a_1 a_3 a_4 a_5+a_1 a_4^2 a_5-(9/2) a_1 a_3 a_5^2-4 a_1 a_4 a_5^2+6 a_1^2 a_4 a_7-(3/2)a_3^3-2 a_3^2 a_4+(9/2) a_3^2 a_5+3 a_3 a_4 a_5-(9/2) a_3 a_5^2-a_4 a_5^2+a_5^3+3 a_1 a_3 a_7+6 a_1 a_4 a_7+3 a_1 a_5 a_7+(7/2) a_3 a_7+a_4 a_7-a_5 a_7
$
\newline\noindent 
$f_5 =
      3 a_1 a_3^3 a_4 a_6+a_1 a_3^2 a_4^2 a_6-9 a_1 a_3^2 a_4 a_5 a_6+8 a_1 a_3 a_4 a_5^2 a_6-2 a_1 a_4^2 a_5^2 a_6-4 a_1 a_4 a_5^3 a_6+5 a_2 a_3^4+a_2 a_3^3 a_4-19 a_2 a_3^3 a_5-2 a_2 a_3^2 a_4 a_5+25 a_2 a_3^2 a_5^2+a_2 a_3 a_4 a_5^2-15 a_2 a_3 a_5^3-a_2 a_4 a_5^3+2 a_2 a_5^4+3 a_3^3 a_4 a_6+a_3^2 a_4^2 a_6-8 a_3^2 a_4 a_5 a_6+7 a_3 a_4 a_5^2 a_6-2 a_4^2 a_5^2 a_6-5 a_4 a_5^3 a_6-9 a_1 a_3 a_4 a_6 a_7+4 a_1 a_4^2 a_6 a_7+9 a_1 a_4 a_5 a_6 a_7-8 a_2 a_3^2 a_7+2 a_2 a_3 a_4 a_7+11 a_2 a_3 a_5 a_7+a_2 a_4 a_5 a_7+11 a_3^2 a_6 a_7-11 a_3 a_4 a_6 a_7+3 a_4^2 a_6 a_7-11 a_3 a_5 a_6 a_7+11 a_4 a_5 a_6 a_7+3 a_5^2 a_6 a_7-2 a_2 a_7^2+2 a_6 a_7^2
 $
\newline\noindent 
$f_6 =
3 a_1^2 a_4^2 a_6+2 a_1 a_3 a_4 a_6+5 a_1 a_4^2 a_6+2 a_1 a_4 a_5 a_6+a_2 a_3^2+2 a_2 a_3 a_4-a_2 a_3 a_5+a_2 a_5^2+2 a_3 a_4 a_6+2 a_4^2 a_6+a_4 a_5 a_6-a_2 a_7+a_6 a_7
$
\newline\noindent 
$f_7 =
      a_1 a_2-a_1 a_6-a_6
$
\newline\noindent 
$f_8 =
     -2 a_1^2 a_2 a_4 a_6^2-a_1 a_2 a_3 a_6^2-3 a_1 a_2 a_4 a_6^2-a_1 a_2 a_5 a_6^2-a_2^2 a_3 a_6-a_2 a_3 a_6^2-a_2 a_4 a_6^2+1
$
\newline\noindent 
$f_9 =
      a_1^2 a_4^3 a_7+a_1 a_3 a_4^2 a_7+a_1 a_4^3 a_7+a_1 a_4^2 a_5 a_7+a_3^2 a_4 a_7+a_3 a_4^2 a_7-a_3 a_4 a_5 a_7+a_4 a_5^2 a_7-a_4 a_7^2+1
 $
\newline\noindent 
$f_{10} =
     5 a_1^2 a_3 a_4^2 a_6-a_1^2 a_4^2 a_5 a_6+3 a_1 a_3^2 a_4 a_6+8 a_1 a_3 a_4^2 a_6+3 a_1 a_3 a_4 a_5 a_6-a_1 a_4^2 a_5 a_6+a_2 a_3^3+3 a_2 a_3^2 a_4+a_2 a_5^3+3 a_3^2 a_4 a_6+3 a_3 a_4^2 a_6+a_3 a_4 a_5 a_6+a_4 a_5^2 a_6-a_2 a_3 a_7-a_2 a_5 a_7
 $
\newline\noindent 
$f_{11} =
     -3 a_1^2 a_4 a_6^2-a_1 a_3 a_6^2-6 a_1 a_4 a_6^2-a_1 a_5 a_6^2-a_2^2 a_3-a_2 a_3 a_6-a_2 a_4 a_6-a_3 a_6^2-3 a_4 a_6^2
 $
\newline\noindent 
$f_{12} =
     a_1^3 a_4^3+2 a_1^2 a_4^3+3 a_1^2 a_4^2 a_5+a_1 a_3 a_4^2+a_1 a_4^3+3 a_1 a_4^2 a_5+3 a_1 a_4 a_5^2+a_3 a_4^2+a_3 a_4 a_5+a_4 a_5^2+a_5^3-2 a_1 a_4 a_7-a_3 a_7-a_4 a_7-a_5 a_7
$
}

\vspace{.2cm}
\normalsize{
  These twelve equations are solvable, and the quantum dimensions, subject to the
  matrix factorisation conditions, are given by 
}

\small{
  $q_L(Q) =  a_1 a_4 a_6+a_2 a_3+a_4 a_6+a_5 a_6$
  
  $q_R(Q) = ( 462  a_1 a_5 a_6^2 a_7^2+603 a_1^3 a_6^2
  -{ 2002} a_2^2 a_3 a_7^2+ {158 } a_2^2 a_4 a_7^2
  -{853} a_2^2 a_5 a_7^2-{ 898 } a_2 a_3 a_6 a_7^2
  - 2784 a_2 a_4 a_6 a_7^2- 136   a_2 a_5 a_6 a_7^2+ 214 a_3 a_6^2 a_7^2
  - 1294 a_4 a_6^2 a_7^2+1111 a_5 a_6^2 a_7^2+ 2646 a_1^2 a_6^2-261 a_1 a_6^2
  -291 a_2^2-301 a_2 a_6-2095 a_6^2)/764 $
}

\normalsize{
  Note that these expressions result after reduction by the ideal spanned
  by the $f_\alpha$, hence the quantum dimensions of this defect are non-zero
  numbers after inserting any special solution to the equations $f_\alpha(a)=0$.
}

\vspace{.5cm}

\normalsize{
\noindent  
    \textbf{(3)} A rank 6 orbifold equivalence between $Z_{13}$ and $Q_{11}$, which
    could be worse:  }

\vspace{.2cm}

\small{

\noindent $ E_{11} =2 y_3 a_1^2 a_2+2 y_3 a_1 a_3 
$ \newline\noindent $ E_{12} =-(3/2)x_1^3 a_1^3 a_2^3-x_1^3 a_1^2 a_3 a_2^2+(1/2) x_1^3 a_1 a_3^2 a_2+2 x_2 y_3 a_4 a_1^2 a_2+2 x_2 y_3 a_4 a_1 a_3-x_2 y_3 a_1+x_1 y_2 a_3+x_3 
$ \newline\noindent $ E_{13} =0  
$ \newline\noindent $ E_{14} =-(3/8) x_1^2 a_4 a_1^3 a_2^3-(1/4) x_1^2 a_4 a_1^2 a_3 a_2^2+(1/8) x_1^2 a_4 a_1 a_3^2 a_2+(1/4) x_1^2 a_3^2-y_2 
$ \newline\noindent $ E_{15} =-x_2 
$ \newline\noindent $ E_{16} =0  
$ \newline\noindent $ E_{21} =(3/2)x_1^3 a_1^3 a_2^3+x_1^3 a_1^2 a_3 a_2^2-(1/2) x_1^3 a_1 a_3^2 a_2+x_2 y_3 a_1-x_1 y_2 a_3+x_3 
$ \newline\noindent $ E_{22} =(3/4) x_1^3 x_2 a_4 a_1^3 a_2^3-x_1^2 y_3^2 a_1^4 a_2^2+(1/2) x_1^3 x_2 a_4 a_1^2 a_3 a_2^2-x_1^2 y_3^2 a_1^3 a_3 a_2-(1/4) x_1^3 x_2 a_4 a_1 a_3^2 a_2-(1/2) x_1^3 x_2 a_3^2+x_2^2 y_3 a_4 a_1+y_2 y_3^2 a_1^2-x_1 x_2 y_2 a_4 a_3+(1/2) x_1^2 y_3^2 a_2-y_1^2 a_5^2+x_1 x_2 y_2+x_2 x_3 a_4 
$ \newline\noindent $ E_{23} =(9/32) x_1^4 a_4^2 a_1^6 a_2^6+(3/8) x_1^4 a_4^2 a_1^5 a_3 a_2^5-(1/16) x_1^4 a_4^2 a_1^4 a_3^2 a_2^4-(1/8) x_1^4 a_4^2 a_1^3 a_3^3 a_2^3+(1/32) x_1^4 a_4^2 a_1^2 a_3^4 a_2^2+(3/8) x_1^2 y_2 a_4 a_1^3 a_2^3+(1/4) x_1^2 y_2 a_4 a_1^2 a_3 a_2^2-(1/8) x_1^4 a_3^4-(1/8) x_1^2 y_2 a_4 a_1 a_3^2 a_2+x_1 x_2 y_3 a_1^2 a_2+(3/4) x_1^2 y_2 a_3^2+x_2 y_1 a_5-y_2^2 
$ \newline\noindent $ E_{24} =0  
$ \newline\noindent $ E_{25} =0  
$ \newline\noindent $ E_{26} =(1/2) x_1^3 y_3 a_1^4 a_2^3-(3/8) x_1^2 y_1 a_4 a_1^3 a_5 a_2^3-(1/4) x_1^2 y_1 a_4 a_1^2 a_3 a_5 a_2^2-(1/2) x_1^3 y_3 a_1^2 a_3^2 a_2+ \newline (1/8) x_1^2 y_1 a_4 a_1 a_3^2 a_5 a_2+(1/4) x_1^2 y_1 a_3^2 a_5+x_1 y_2 y_3 a_1^2 a_2+x_2 y_3^2 a_1^2+x_1 x_2^2-y_1 y_2 a_5 
$ \newline\noindent $ E_{31} =0  
$ \newline\noindent $ E_{32} =(9/32) x_1^4 a_4^2 a_1^6 a_2^6+(3/8) x_1^4 a_4^2 a_1^5 a_3 a_2^5-(1/16) x_1^4 a_4^2 a_1^4 a_3^2 a_2^4-(1/8) x_1^4 a_4^2 a_1^3 a_3^3 a_2^3+(1/32) x_1^4 a_4^2 a_1^2 a_3^4 a_2^2+(3/8) x_1^2 y_2 a_4 a_1^3 a_2^3+(1/4) x_1^2 y_2 a_4 a_1^2 a_3 a_2^2-(1/8) x_1^4 a_3^4-(1/8) x_1^2 y_2 a_4 a_1 a_3^2 a_2+x_1 x_2 y_3 a_1^2 a_2+(3/4) x_1^2 y_2 a_3^2-x_2 y_1 a_5-y_2^2 
$ \newline\noindent $ E_{33} =-(3/4) x_1^2 y_3 a_4 a_1^5 a_2^4-(1/2) x_1^2 y_3 a_4 a_1^4 a_3 a_2^3+(1/4) x_1^2 y_3 a_4 a_1^3 a_3^2 a_2^2+3 x_1^2 y_3 a_1^4 a_2^3+2 x_1^2 y_3 a_1^3 a_3 a_2^2-(1/2) x_1^2 y_3 a_1^2 a_3^2 a_2-2 y_2 y_3 a_1^2 a_2-2 y_2 y_3 a_1 a_3+x_2^2 
$ \newline\noindent $ E_{34} =(3/2)x_1^3 a_1^3 a_2^3+x_1^3 a_1^2 a_3 a_2^2-(1/2) x_1^3 a_1 a_3^2 a_2+x_2 y_3 a_1-x_1 y_2 a_3+x_3 
$ \newline\noindent $ E_{35} =0  
$ \newline\noindent $ E_{36} =-(3/8) x_1^2 x_2 a_4 a_1^3 a_2^3+2 x_1 y_3^2 a_1^4 a_2^2-(1/4) x_1^2 x_2 a_4 a_1^2 a_3 a_2^2+2 x_1 y_3^2 a_1^3 a_3 a_2+(1/8) x_1^2 x_2 a_4 a_1 a_3^2 a_2-2 y_1 y_3 a_1^2 a_5 a_2+(1/4) x_1^2 x_2 a_3^2-2 y_1 y_3 a_1 a_3 a_5-x_2 y_2 
$ \newline\noindent $ E_{41} =-(3/8) x_1^2 a_4 a_1^3 a_2^3-(1/4) x_1^2 a_4 a_1^2 a_3 a_2^2+(1/8) x_1^2 a_4 a_1 a_3^2 a_2+(1/4) x_1^2 a_3^2-y_2 
$ \newline\noindent $ E_{42} =-(3/8) x_1^2 x_2 a_4^2 a_1^3 a_2^3-(1/4) x_1^2 x_2 a_4^2 a_1^2 a_3 a_2^2+(1/8) x_1^2 x_2 a_4^2 a_1 a_3^2 a_2+(1/4) x_1^2 x_2 a_4 a_3^2-x_2 y_2 a_4 
$ \newline\noindent $ E_{43} =-(3/2)x_1^3 a_1^3 a_2^3-x_1^3 a_1^2 a_3 a_2^2+(1/2) x_1^3 a_1 a_3^2 a_2-x_2 y_3 a_1+x_1 y_2 a_3+x_3 
$ \newline\noindent $ E_{44} =-y_3^2 a_1^2-x_1 x_2 
$ \newline\noindent $ E_{45} =x_1 y_3 a_1^2 a_2-y_1 a_5 
$ \newline\noindent $ E_{46} =0  
$ \newline\noindent $ E_{51} =-x_2 
$ \newline\noindent $ E_{52} =-x_2^2 a_4 
$ \newline\noindent $ E_{53} =0  
$ \newline\noindent $ E_{54} =x_1 y_3 a_1^2 a_2+y_1 a_5 
$ \newline\noindent $ E_{55} =-x_1^2 a_1^2 a_2^2+y_2 
$ \newline\noindent $ E_{56} =-(3/2)x_1^3 a_1^3 a_2^3-x_1^3 a_1^2 a_3 a_2^2+(1/2) x_1^3 a_1 a_3^2 a_2-x_2 y_3 a_1+x_1 y_2 a_3+x_3 
$ \newline\noindent $ E_{61} =-x_1 x_2 a_3 
$ \newline\noindent $ E_{62} =(1/2) x_1^3 y_3 a_1^4 a_2^3+(3/8) x_1^2 y_1 a_4 a_1^3 a_5 a_2^3+(1/4) x_1^2 y_1 a_4 a_1^2 a_3 a_5 a_2^2-(1/2) x_1^3 y_3 a_1^2 a_3^2 a_2- \newline (1/8) x_1^2 y_1 a_4 a_1 a_3^2 a_5 a_2-(1/4) x_1^2 y_1 a_3^2 a_5+x_1 y_2 y_3 a_1^2 a_2+x_2 y_3^2 a_1^2-x_1 x_2^2 a_4 a_3+x_1 x_2^2+y_1 y_2 a_5 
$ \newline\noindent $ E_{63} =-(3/8) x_1^2 x_2 a_4 a_1^3 a_2^3-(1/4) x_1^2 x_2 a_4 a_1^2 a_3 a_2^2+(1/8) x_1^2 x_2 a_4 a_1 a_3^2 a_2+2 y_1 y_3 a_1^2 a_5 a_2+(1/4) x_1^2 x_2 a_3^2+2 y_1 y_3 a_1 a_3 a_5+x_1 y_3^2 a_2-x_2 y_2 
$ \newline\noindent $ E_{64} =x_1^2 y_3 a_1^2 a_3 a_2+x_1 y_1 a_3 a_5 
$ \newline\noindent $ E_{65} =(3/2)x_1^3 a_1^3 a_2^3-(1/2) x_1^3 a_1 a_3^2 a_2+x_2 y_3 a_1+x_3 
$ \newline\noindent $ E_{66} =x_1^4 a_1^4 a_2^4+(3/2)x_1^4 a_1^3 a_3 a_2^3-(1/2) x_1^4 a_1 a_3^3 a_2+2 y_3^3 a_1^4 a_2+2 y_3^3 a_1^3 a_3+x_1^2 y_2 a_1^2 a_2^2+2 x_1 x_2 y_3 a_1^2 a_2+x_1 x_2 y_3 a_1 a_3+x_1 x_3 a_3+y_2^2 
$

\vspace{.2cm}

\noindent $ J_{11} =-(3/4) x_1^3 x_2 a_4 a_1^3 a_2^3+x_1^2 y_3^2 a_1^4 a_2^2-(1/2) x_1^3 x_2 a_4 a_1^2 a_3 a_2^2+x_1^2 y_3^2 a_1^3 a_3 a_2+(1/4) x_1^3 x_2 a_4 a_1 a_3^2 a_2+(1/2) x_1^3 x_2 a_3^2-x_2^2 y_3 a_4 a_1-y_2 y_3^2 a_1^2+x_1 x_2 y_2 a_4 a_3-(1/2) x_1^2 y_3^2 a_2+y_1^2 a_5^2-x_1 x_2 y_2-x_2 x_3 a_4 
$ \newline\noindent $ J_{12} =-(3/2)x_1^3 a_1^3 a_2^3-x_1^3 a_1^2 a_3 a_2^2+(1/2) x_1^3 a_1 a_3^2 a_2+2 x_2 y_3 a_4 a_1^2 a_2+2 x_2 y_3 a_4 a_1 a_3-x_2 y_3 a_1+x_1 y_2 a_3+x_3 
$ \newline\noindent $ J_{13} =-(3/8) x_1^2 x_2 a_4^2 a_1^3 a_2^3-(1/4) x_1^2 x_2 a_4^2 a_1^2 a_3 a_2^2+(1/8) x_1^2 x_2 a_4^2 a_1 a_3^2 a_2+(1/4) x_1^2 x_2 a_4 a_3^2-x_2 y_2 a_4 
$ \newline\noindent $ J_{14} =-(9/32) x_1^4 a_4^2 a_1^6 a_2^6-(3/8) x_1^4 a_4^2 a_1^5 a_3 a_2^5+(1/16) x_1^4 a_4^2 a_1^4 a_3^2 a_2^4+(1/8) x_1^4 a_4^2 a_1^3 a_3^3 a_2^3-(1/32) x_1^4 a_4^2 a_1^2 a_3^4 a_2^2-(3/8) x_1^2 y_2 a_4 a_1^3 a_2^3-(1/4) x_1^2 y_2 a_4 a_1^2 a_3 a_2^2+(1/8) x_1^4 a_3^4+(1/8) x_1^2 y_2 a_4 a_1 a_3^2 a_2-x_1 x_2 y_3 a_1^2 a_2-(3/4) x_1^2 y_2 a_3^2-x_2 y_1 a_5+y_2^2 
$ \newline\noindent $ J_{15} =-(1/2) x_1^3 y_3 a_1^4 a_2^3+(3/8) x_1^2 y_1 a_4 a_1^3 a_5 a_2^3+(1/4) x_1^2 y_1 a_4 a_1^2 a_3 a_5 a_2^2+(1/2) x_1^3 y_3 a_1^2 a_3^2 a_2-(1/8) x_1^2 y_1 a_4 a_1 a_3^2 a_5 a_2-(1/4) x_1^2 y_1 a_3^2 a_5-x_1 y_2 y_3 a_1^2 a_2-x_2 y_3^2 a_1^2+x_1 x_2^2 a_4 a_3-x_1 x_2^2+y_1 y_2 a_5 
$ \newline\noindent $ J_{16} =-x_2^2 a_4 
$ \newline\noindent $ J_{21} =(3/2)x_1^3 a_1^3 a_2^3+x_1^3 a_1^2 a_3 a_2^2-(1/2) x_1^3 a_1 a_3^2 a_2+x_2 y_3 a_1-x_1 y_2 a_3+x_3 
$ \newline\noindent $ J_{22} =-2 y_3 a_1^2 a_2-2 y_3 a_1 a_3 
$ \newline\noindent $ J_{23} =(3/8) x_1^2 a_4 a_1^3 a_2^3+(1/4) x_1^2 a_4 a_1^2 a_3 a_2^2-(1/8) x_1^2 a_4 a_1 a_3^2 a_2-(1/4) x_1^2 a_3^2+y_2 
$ \newline\noindent $ J_{24} =0  
$ \newline\noindent $ J_{25} =-x_1 x_2 a_3 
$ \newline\noindent $ J_{26} =x_2 
$ \newline\noindent $ J_{31} =0  
$ \newline\noindent $ J_{32} =(3/8) x_1^2 a_4 a_1^3 a_2^3+(1/4) x_1^2 a_4 a_1^2 a_3 a_2^2-(1/8) x_1^2 a_4 a_1 a_3^2 a_2-(1/4) x_1^2 a_3^2+y_2 
$ \newline\noindent $ J_{33} =y_3^2 a_1^2+x_1 x_2 
$ \newline\noindent $ J_{34} =(3/2)x_1^3 a_1^3 a_2^3+x_1^3 a_1^2 a_3 a_2^2-(1/2) x_1^3 a_1 a_3^2 a_2+x_2 y_3 a_1-x_1 y_2 a_3+x_3 
$ \newline\noindent $ J_{35} =x_1^2 y_3 a_1^2 a_3 a_2-x_1 y_1 a_3 a_5 
$ \newline\noindent $ J_{36} =-x_1 y_3 a_1^2 a_2+y_1 a_5 
$ \newline\noindent $ J_{41} =-(9/32) x_1^4 a_4^2 a_1^6 a_2^6-(3/8) x_1^4 a_4^2 a_1^5 a_3 a_2^5+(1/16) x_1^4 a_4^2 a_1^4 a_3^2 a_2^4+(1/8) x_1^4 a_4^2 a_1^3 a_3^3 a_2^3-   \newline (1/32) x_1^4 a_4^2 a_1^2 a_3^4 a_2^2-(3/8) x_1^2 y_2 a_4 a_1^3 a_2^3-(1/4) x_1^2 y_2 a_4 a_1^2 a_3 a_2^2+(1/8) x_1^4 a_3^4+ \newline (1/8) x_1^2 y_2 a_4 a_1 a_3^2 a_2-x_1 x_2 y_3 a_1^2 a_2-(3/4) x_1^2 y_2 a_3^2+x_2 y_1 a_5+y_2^2 
$ \newline\noindent $ J_{42} =0  
$ \newline\noindent $ J_{43} =-(3/2)x_1^3 a_1^3 a_2^3-x_1^3 a_1^2 a_3 a_2^2+(1/2) x_1^3 a_1 a_3^2 a_2-x_2 y_3 a_1+x_1 y_2 a_3+x_3 
$ \newline\noindent $ J_{44} =(3/4) x_1^2 y_3 a_4 a_1^5 a_2^4+(1/2) x_1^2 y_3 a_4 a_1^4 a_3 a_2^3-(1/4) x_1^2 y_3 a_4 a_1^3 a_3^2 a_2^2-3 x_1^2 y_3 a_1^4 a_2^3-2 x_1^2 y_3 a_1^3 a_3 a_2^2+(1/2) x_1^2 y_3 a_1^2 a_3^2 a_2+2 y_2 y_3 a_1^2 a_2+2 y_2 y_3 a_1 a_3-x_2^2 
$ \newline\noindent $ J_{45} =(3/8) x_1^2 x_2 a_4 a_1^3 a_2^3-2 x_1 y_3^2 a_1^4 a_2^2+(1/4) x_1^2 x_2 a_4 a_1^2 a_3 a_2^2-2 x_1 y_3^2 a_1^3 a_3 a_2-(1/8) x_1^2 x_2 a_4 a_1 a_3^2 a_2+2 y_1 y_3 a_1^2 a_5 a_2-(1/4) x_1^2 x_2 a_3^2+2 y_1 y_3 a_1 a_3 a_5+x_2 y_2 
$ \newline\noindent $ J_{46} =0  
$ \newline\noindent $ J_{51} =-(1/2) x_1^3 y_3 a_1^4 a_2^3-(3/8) x_1^2 y_1 a_4 a_1^3 a_5 a_2^3-(1/4) x_1^2 y_1 a_4 a_1^2 a_3 a_5 a_2^2+(1/2) x_1^3 y_3 a_1^2 a_3^2 a_2+(1/8) x_1^2 y_1 a_4 a_1 a_3^2 a_5 a_2+(1/4) x_1^2 y_1 a_3^2 a_5-x_1 y_2 y_3 a_1^2 a_2-x_2 y_3^2 a_1^2-x_1 x_2^2-y_1 y_2 a_5 
$ \newline\noindent $ J_{52} =0  
$ \newline\noindent $ J_{53} =0  
$ \newline\noindent $ J_{54} =(3/8) x_1^2 x_2 a_4 a_1^3 a_2^3+(1/4) x_1^2 x_2 a_4 a_1^2 a_3 a_2^2-(1/8) x_1^2 x_2 a_4 a_1 a_3^2 a_2-2 y_1 y_3 a_1^2 a_5 a_2-(1/4) x_1^2 x_2 a_3^2-2 y_1 y_3 a_1 a_3 a_5-x_1 y_3^2 a_2+x_2 y_2 
$ \newline\noindent $ J_{55} =-x_1^4 a_1^4 a_2^4-(3/2)x_1^4 a_1^3 a_3 a_2^3+(1/2) x_1^4 a_1 a_3^3 a_2-2 y_3^3 a_1^4 a_2-2 y_3^3 a_1^3 a_3-x_1^2 y_2 a_1^2 a_2^2-2 x_1 x_2 y_3 a_1^2 a_2-x_1 x_2 y_3 a_1 a_3-x_1 x_3 a_3-y_2^2 
$ \newline\noindent $ J_{56} =-(3/2)x_1^3 a_1^3 a_2^3-x_1^3 a_1^2 a_3 a_2^2+(1/2) x_1^3 a_1 a_3^2 a_2-x_2 y_3 a_1+x_1 y_2 a_3+x_3 
$ \newline\noindent $ J_{61} =0  
$ \newline\noindent $ J_{62} =x_2 
$ \newline\noindent $ J_{63} =-x_1 y_3 a_1^2 a_2-y_1 a_5 
$ \newline\noindent $ J_{64} =0  
$ \newline\noindent $ J_{65} =(3/2)x_1^3 a_1^3 a_2^3-(1/2) x_1^3 a_1 a_3^2 a_2+x_2 y_3 a_1+x_3 
$ \newline\noindent $ J_{66} =x_1^2 a_1^2 a_2^2-y_2 
$
}

\vspace{.2cm}
\normalsize{ 
  The five coefficients $a_1,\ldots,a_5$ are subject to thirty-seven
  relatively simple conditions $f_\alpha(a)=0$ with 
}

\vspace{.2cm}
\small{
\noindent $f_{1} = a_1^2+a_5^2 
$ \newline\noindent $f_{2} = -3 a_1 a_2 a_4+52 a_1 a_3 a_5^2-7 a_3 a_4-10 
$ \newline\noindent $f_{3} = -1839 a_1 a_2 a_3+30 a_2^2 a_5^2+835 a_3^2+72 a_4^4 
$ \newline\noindent $f_{4} = 94888 a_1 a_4 a_5^2-6675 a_2^3+7504 a_3 a_4^3+41908 a_4^2 
$ \newline\noindent $f_{5} = -159 a_1 a_2 a_4+52 a_3^2 a_4^2+383 a_3 a_4+445 
$ \newline\noindent $f_{6} = 83 a_1 a_2 a_5^2+14 a_2 a_3 a_4^2+53 a_2 a_4-75 a_3 a_5^2 
$ \newline\noindent $f_{7} = 225 a_1 a_2 a_4^2+2314 a_1 a_5^2+187 a_3 a_4^2+724 a_4 
$ \newline\noindent $f_{8} = -36 a_1 a_2 a_3+15 a_2^2 a_5^2+8 a_3^3 a_4+47 a_3^2 
$ \newline\noindent $f_{9} = 2 a_1 a_3^2 a_4+78 a_1 a_3-15 a_2^2 a_4^2-81 a_2 a_5^2 
$ \newline\noindent $f_{10} = 10 a_1 a_2 a_3 a_4+13 a_1 a_2+6 a_3^2 a_4+29 a_3 
$ \newline\noindent $f_{11} = 33 a_1 a_2^2 a_4-27 a_2 a_3 a_4-72 a_2-52 a_3^2 a_5^2 
$ \newline\noindent $f_{12} = -145863 a_2^3+47444 a_3^4-31896 a_3 a_4^3-82080 a_4^2 
$ \newline\noindent $f_{13} = 1892 a_1 a_3^3-648 a_1 a_4^3+228 a_2^2 a_3 a_4-2223 a_2^2 
$ \newline\noindent $f_{14} = 8 a_1 a_3 a_4^2+44 a_1 a_4+a_2^2 a_3^2-92 a_5^4 
$ \newline\noindent $f_{15} = 1311 a_1 a_2 a_3^2+342 a_2^3 a_4-211 a_3^3+180 a_4^3 
$ \newline\noindent $f_{16} = -4926 a_1 a_2 a_4^2+1157 a_2^3 a_3+102 a_3 a_4^2+7968 a_4 
$ \newline\noindent $f_{17} = 41 a_2^4-132 a_2 a_4^2-408 a_3 a_4 a_5^2-784 a_5^2 
$ \newline\noindent $f_{18} = 1002 a_1 a_2^3-12 a_1 a_4^2-469 a_2^2 a_3+1224 a_4 a_5^4 
$ \newline\noindent $f_{19} = 7 a_1 a_3 a_4+10 a_1-3 a_2 a_4 a_5^2+52 a_3 a_5^4 
$ \newline\noindent $f_{20} = -2 a_1 a_3 a_5^2+2 a_2 a_5^4-1 
$ \newline\noindent $f_{21} = 67716 a_1 a_5^4-959 a_2 a_3^3-1584 a_2 a_4^3+23256 a_4 a_5^2 
$ \newline\noindent $f_{22} = -71 a_1 a_2^2+116 a_2 a_3^2 a_4+403 a_2 a_3+48 a_4^3 a_5^2 
$ \newline\noindent $f_{23} = 3649 a_2 a_3^3+9999 a_2 a_4^3+33858 a_3 a_4^2 a_5^2+55575 a_4 a_5^2 
$ \newline\noindent $f_{24} = 654 a_1 a_3 a_4^2+17604 a_1 a_4+1157 a_2^2 a_3^2+10350 a_2 a_4^2 a_5^2 
$ \newline\noindent $f_{25} = 283176 a_1 a_4^2 a_5^2-170487 a_2^3 a_4-166964 a_3^3+46476 a_4^3 
$ \newline\noindent $f_{26} = 123 a_1 a_2 a_5^2+15 a_2 a_3 a_4^2+106 a_3^2 a_4 a_5^2+514 a_3 a_5^2 
$ \newline\noindent $f_{27} = -6 a_1 a_3^2 a_4-29 a_1 a_3+10 a_2 a_3 a_4 a_5^2+13 a_2 a_5^2 
$ \newline\noindent $f_{28} = 3 a_1 a_2+3 a_2^2 a_4 a_5^2-a_3^2 a_4-7 a_3 
$ \newline\noindent $f_{29} = 246 a_1 a_2 a_4 a_5^2-39 a_2 a_4^2-128 a_3 a_4 a_5^2-623 a_5^2 
$ \newline\noindent $f_{30} = 33 a_1 a_2^2+36 a_2 a_3^2 a_4+129 a_2 a_3+40 a_3^3 a_5^2 
$ \newline\noindent $f_{31} = -22 a_1 a_3^3-24 a_2^2 a_3 a_4-9 a_2^2+54 a_2 a_3^2 a_5^2 
$ \newline\noindent $f_{32} = 2 a_1 a_2 a_3^2-6 a_2^3 a_4+45 a_2^2 a_3 a_5^2+3 a_3^3 
$ \newline\noindent $f_{33} = -40 a_1 a_2^2 a_3+123 a_2^3 a_5^2+51 a_2 a_3^2+48 a_4^2 a_5^2 
$ \newline\noindent $f_{34} = 9 a_1 a_2^2 a_5^2+7 a_1 a_3^2+6 a_2^2 a_4-18 a_2 a_3 a_5^2 
$ \newline\noindent $f_{35} = -329 a_1 a_3 a_4+89 a_1+78 a_2^2 a_4^3+453 a_2 a_4 a_5^2 
$ \newline\noindent $f_{36} = -3211 a_1 a_2 a_3+393 a_2^3 a_4^2+2637 a_2^2 a_5^2+306 a_3^3 a_4+1765 a_3^2 
$ \newline\noindent $f_{37} = -5112 a_1 a_4 a_5^2+19 a_2^3+612 a_4^2+15008 a_5^6 
$
}

\vspace{.2cm}
\normalsize{
  The quantum dimensions of $Q$ are:
  
 $q_L(Q) = (24/13) a_2 a_4 a_5^3-(4/13) a_1 a_3 a_4 a_5+(50/13) a_1 a_5$ 

 $q_R(Q) = -2 a_1 a_2 a_5-2 a_3 a_5$ . 
}

\vspace{.5cm}

\normalsize{
\noindent
\textbf{(4)} A rank 4 orbifold equivalence between $S_{11}$ and $W_{13}$:  }

\vspace{.2cm}   

\small{
\noindent $ E_{11} =x_1 y_3-y_1 $
\newline\noindent $ E_{12} =-x_2 y_3^2+x_1^2+x_2 x_3 $
\newline\noindent $ E_{13} =-x_2+y_2 $
\newline\noindent $ E_{14} =0 $
\newline\noindent $ E_{21} =-y_3^2-x_3 $
\newline\noindent $ E_{22} =-x_1 y_3-y_1 $
\newline\noindent $ E_{23} =0 $
\newline\noindent $ E_{24} =-x_2+y_2 $
\newline\noindent $ E_{31} =x_2^3+x_2^2 y_2+x_2 y_2^2+y_2^3-x_3 y_3^2 $
\newline\noindent $ E_{32} =-x_1 y_3^3-y_1 y_3^2 $
\newline\noindent $ E_{33} =-x_1 y_3-y_1 $
\newline\noindent $ E_{34} =y_2 y_3^2-x_1^2-x_2 x_3 $
\newline\noindent $ E_{41} =0 $
\newline\noindent $ E_{42} =y_3^4+x_2^3+x_2^2 y_2+x_2 y_2^2+y_2^3 $
\newline\noindent $ E_{43} =y_3^2+x_3 $
\newline\noindent $ E_{44} =x_1 y_3-y_1 $

\vspace{.2cm}

\noindent $ J_{11} =-x_1 y_3-y_1 $
\newline\noindent $ J_{12} =y_2 y_3^2-x_1^2-x_2 x_3 $
\newline\noindent $ J_{13} =x_2-y_2 $
\newline\noindent $ J_{14} =0 $
\newline\noindent $ J_{21} =y_3^2+x_3 $
\newline\noindent $ J_{22} =x_1 y_3-y_1 $
\newline\noindent $ J_{23} =0 $
\newline\noindent $ J_{24} =x_2-y_2 $
\newline\noindent $ J_{31} =-y_3^4-x_2^3-x_2^2 y_2-x_2 y_2^2-y_2^3 $
\newline\noindent $ J_{32} =-x_1 y_3^3+y_1 y_3^2 $
\newline\noindent $ J_{33} =x_1 y_3-y_1 $
\newline\noindent $ J_{34} =-x_2 y_3^2+x_1^2+x_2 x_3 $
\newline\noindent $ J_{41} =0 $
\newline\noindent $ J_{42} =-x_2^3-x_2^2 y_2-x_2 y_2^2-y_2^3+x_3 y_3^2 $
\newline\noindent $ J_{43} =-y_3^2-x_3 $
\newline\noindent $ J_{44} =-x_1 y_3-y_1 $
}

\vspace{.2cm} 

\normalsize{ 
 This rather simple $Q$ does not depend on any coefficients, although more general
 orbifold equivalences between $S_{11}$ and $W_{13}$ can be found. 

 Its quantum dimensions are $q_L(Q) = -2$ and $q_R(Q) =  -1$.
}

\vspace{.5cm}

\noindent  \normalsize{
\textbf{(5)} A rank 3 orbifold equivalence between a chain and a loop at
   central charge $\hat c = \frac{6}{5}\,$: }

\vspace{.2cm}   
\small{
\noindent $ E_{11}= 2 a_1  x_1^4+ 2 a_1  x_1 y_2^2 $
\newline\noindent $ E_{12}= a_1  x_1^3 y_2+ a_1  y_2^3+y_1 $
\newline\noindent $ E_{13}= a_1  x_1^5+ a_1  x_1^2 y_2^2+x_2 $
\newline\noindent $ E_{21}= -2 a_1  x_1^3 y_2-a_1  y_2^3+y_1 $
\newline\noindent $ E_{22}= -a_1  x_1^2 y_2^2+x_2 $
\newline\noindent $ E_{23}= -a_1  x_1^4 y_2+x_1 y_1 $
\newline\noindent $ E_{31}= -a_1  x_1^5+x_2 $
\newline\noindent $ E_{32}=-x_1 y_1 $
\newline\noindent $ E_{33}=-y_1 y_2 $

\vspace{.2cm}
\noindent $ J_{11}= a_1  x_1^5 y_1 y_2 -a_1  x_1^2 y_1 y_2^3-x_1^2 y_1^2+x_2 y_1 y_2 $
\newline\noindent $ J_{12}= a_1  x_1^6 y_1-a_1  y_1 y_2^4+x_1 x_2 y_1-y_1^2 y_2 $
\newline\noindent $ J_{13}= a_1^2+1  x_1^{10}+ a_1  x_1^5 x_2-a_1  x_1 y_1 y_2^3-x_1 y_1^2+x_2^2 $
\newline\noindent $ J_{21}= -a_1^2  x_1^9 y_2+ a_1  x_1^6 y_1+ a_1  x_1^4 x_2 y_2+ 2 a_1  x_1^3 y_1 y_2^2+ a_1  y_1 y_2^4-x_1 x_2 y_1-y_1^2 y_2 $
\newline\noindent $ J_{22}=x_1^{10}-a_1^2  x_1^7 y_2^2+ 2 a_1  x_1^4 y_1 y_2+ a_1  x_1^2 x_2 y_2^2+ 2 a_1  x_1 y_1 y_2^3+x_2^2 $
\newline\noindent $ J_{23}= a_1^2  x_1^5 y_2^3+ a_1^2  x_1^2 y_2^5+ a_1  x_1^5 y_1+ 2 a_1  x_1^3 x_2 y_2+ a_1  x_1^2 y_1 y_2^2+ a_1  x_2 y_2^3-x_2 y_1 $
\newline\noindent $ J_{31}= (a_1^2+1)  x_1^{10}+ a_1^2  x_1^7 y_2^2-a_1  x_1^5 x_2-2 a_1  x_1^4 y_1 y_2-a_1  x_1^2 x_2 y_2^2-a_1  x_1 y_1 y_2^3+x_1 y_1^2+x_2^2 $
\newline\noindent $ J_{32}= a_1^2  x_1^8 y_2+ a_1^2  x_1^5 y_2^3-a_1  x_1^5 y_1-a_1  x_1^3 x_2 y_2-2 a_1  x_1^2 y_1 y_2^2-a_1  x_2 y_2^3-x_2 y_1 $
\newline\noindent $ J_{33}= -a_1^2  x_1^3 y_2^4-a_1^2  y_2^6-2 a_1  x_1^4 x_2-a_1  x_1^3 y_1 y_2-2 a_1  x_1 x_2 y_2^2+y_1^2 $

}

\normalsize{
\vspace{.2cm}
This contains only a single coefficient $a_1$ which has to satisfy $a_1^2 = -1$.

\vspace{.2cm}
The quantum dimensions of this defect are $q_L(Q)=-2$ and $q_R(Q)=-3\,$.
}

\end{appendix}

\newpage
\normalsize{

}

\end{document}